\documentclass[12pt,a4paper]{article}
\usepackage{latexsym,epsfig,fullpage,amsfonts,amssymb,amsmath,amscd,epic,graphics,rotating,theorem}
\theoremstyle{plain}
\newtheorem{lemma}{Lemma}
\newtheorem{proposition}[lemma]{Proposition}
\newtheorem{theorem}{Theorem}
\newtheorem{conjecture}{Conjecture}

{\theorembodyfont{\rmfamily}
\newtheorem{definition}{Definition}
\newtheorem{remark}{Remark}
}
\begin{document}
\newcommand{\pperp}{\hbox{$\perp\hskip-6pt\perp$}}
\newcommand{\N}{{\mathbb N}}
\newcommand{\PP}{{\mathbb P}}
\newcommand{\Z}{{\mathbb Z}}
\newcommand{\R}{{\mathbb R}}
\newcommand{\C}{{\mathbb C}}
\newcommand{\Q}{{\mathbb Q}}
\newcommand{\proofend}{$\Box$\bigskip}
\newcommand{\eps}{{\varepsilon}}
\newcommand{\ko}{{\cal O}}
\newcommand{\mt}{{\operatorname{mt}}}
\newcommand{\Def}{{\operatorname{Def}}}
\newcommand{\oT}{{\overline T}}
\newcommand{\ord}{{\operatorname{ord}}}
\newcommand{\Tran}{{\operatorname{Tran}}}
\newcommand{\jet}{{\operatorname{jet}}}
\newcommand{\Iso}{{\operatorname{Iso}}}
\newcommand{\Ker}{{\operatorname{Ker}}}
\newcommand{\bi}{{\bf i}}
\newcommand{\bj}{{\bf j}}
\newcommand{\bk}{{\bf k}}
\newcommand{\bm}{{\bf m}}
\newcommand{\bx}{{\bf x}}
\newcommand{\const}{{\operatorname{const}}}
\newcommand{\conv}{{\operatorname{conv}}}
\newcommand{\Sing}{{\operatorname{Sing}}}
\newcommand{\conj}{{\operatorname{Conj}}}
\newcommand{\Aut}{{\operatorname{Aut}}}
\newcommand{\Sym}{{\operatorname{Sym}}}
\newcommand{\ow}{{\overline w}}
\newcommand{\ov}{{\overline v}}
\newcommand{\kf}{{\cal F}}
\newcommand{\kc}{{\cal C}}
\newcommand{\ki}{{\cal I}}
\newcommand{\kj}{{\cal J}}
\newcommand{\ke}{{\cal E}}
\newcommand{\kz}{{\cal Z}}
\newcommand{\tet}{{\theta}}
\newcommand{\Del}{{\Delta}}
\newcommand{\bet}{{\beta}}
\newcommand{\mm}{{\mathfrak m}}
\newcommand{\kap}{{\kappa}}
\newcommand{\del}{{\delta}}
\newcommand{\sig}{{\sigma}}
\newcommand{\alp}{{\alpha}}
\newcommand{\Sig}{{\Sigma}}
\newcommand{\Gam}{{\Gamma}}
\newcommand{\gam}{{\gamma}}
\newcommand{\Lam}{{\Lambda}}
\newcommand{\lam}{{\lambda}}
\title{Analytic order
of singular and critical points}
\author{E. Shustin
\thanks{The
author was partially supported by Grant No. G-616-15.6/99 of the
German-Israeli Foundation for Research and Development and by the
Hermann-Minkowski Minerva Center for Geometry at Tel Aviv
University. This work has been completed during the author's RiP stay
at the Mathematisches Forschunsinstitut Oberwolfach.}
\\
School of Mathematical Sciences\\
Tel Aviv University\\
Ramat Aviv, 69978 Tel Aviv, Israel\\
E-mail: shustin@post.tau.ac.il}
\date{}
\maketitle
\begin{abstract}
We deal with the following closely related problems: (i) For a
germ of a reduced plane analytic curve, what is the minimal degree
of an algebraic curve with a singular point analytically
equivalent (isomorphic) to the given one? (ii) For a germ of a
holomorphic function in two variables with an isolated critical
point, what is the minimal degree of a polynomial, equivalent to
the given function up to a local holomorphic coordinate change?
Classically known estimates for such a degree $d$ in these
questions are $\sqrt{\mu}+1\le d\le \mu+1$, where $\mu$ is the
Milnor number. Our result in both the problems is $d\le
a\sqrt{\mu}$ with an absolute constant $a$. As a corollary, we
obtain asymptotically proper sufficient conditions for the
existence of algebraic curves with prescribed singularities on
smooth algebraic surfaces.
\end{abstract}

\section*{Introduction}

We work with algebraic curves, hypersurfaces, functions over the
complex field $\C$, though main results hold for any algebraically
closed field of characteristic $0$ by the Lefschetz principle.

{\bf Statement of the problem and formulation of main results}. A
classical question is how many and what types of singular points
can occur on a plane algebraic curve of a given degree. The
well-known bound
$$n\le\frac{(d-1)(d-2)}{2}$$
is necessary and sufficient for the existence of an irreducible
plane curve of degree $d$ with $n$ nodes \cite{Se}. Even the case
of cusps appears to be much more difficult (see \cite{Hi,HI,Sh1}),
so that one cannot expect a complete answer. However, we can ask
for a reasonable sufficient existence conditions which cover
arbitrary degrees and singularities. Namely, the inequalities
$$\mu_0\le(d-1)^2,\quad\del_0\le\frac{(d-1)(d-2)}{2}$$
are necessary for the existence of an irreducible plane curve of
degree $d$ with given singularities having the total Milnor number
$\mu_0$ and the total $\del$-invariant $\del_0$. The main result
of \cite{GLS} (refined later in \cite{Lo}) states that the
condition
\begin{equation}
\mu_0\le\frac{1}{46}(d+2)^2\label{e90}
\end{equation}
is sufficient for the existence of an irreducible plane curve of degree $d$
with given singularities prescribed up to {\it topological}
equivalence. Asymptotically it coincides with the necessary condition
up to a constant factor, and, thus, is called {\it asymptotically proper}.

However, (\ref{e90}) does not apply to singularities defined up to
{\it analytic} equivalence\footnote{Also called {\it analytic
isomorphism} or {\it contact} equivalence.}. Our first result
(Theorem \ref{t3} and Remark \ref{r2}, section \ref{sec6}) is that
the inequality
$$\mu_0\le\frac{1}{9}(d^2-2d+3)$$
is sufficient for the existence of an irreducible plane curve of
degree $d$ with arbitrary singularities prescribed up to {\it
analytic} equivalence and with the total Milnor number $\mu_0$.
This not only covers a wider range than (\ref{e90}), but is
concerned with a stronger equivalence relation for singular
points. In the case of one singularity represented by a curve germ
$(C,z)$ we estimate the {\it analytic order} $d(C,z)$ of this
germ, i.e., the minimal degree of a plane curve having a singular
point analytically equivalent to $(C,z)$, as
$$d(C,z)\le 3\sqrt{\mu(C,z)}-1$$
(Theorem \ref{t2} and Remark \ref{r2}, section \ref{sec6}). A
closely related question: given a germ of a holomorphic function
$f:(\C^2,0)\to(\C,0)$ with an isolated critical point, what is the
{\it analytic order} $d(f)$ of this germ, i.e., the minimal degree
of a polynomial equivalent to $f$ up to a local holomorphic
coordinate change (so-called {\it right} equivalence)? We refine
the classical bounds
$$\sqrt{\mu(f)}+1\le d(f)\le\mu(f)+1$$
with our upper bound
$$d(f)<4\sqrt{\mu(f)}-1$$
(Theorem \ref{t5}, section \ref{sec3}).

These questions can be generalized in two directions. First, one
can look for curves with prescribed singularities in given linear
systems on smooth algebraic surfaces. We provide a sufficient
numerical condition for the existence of an irreducible curve with
singularities prescribed up to analytic equivalence in a given
linear system on a smooth algebraic surface (Theorem \ref{t4},
section \ref{sec4}). It is stronger than a similar sufficient
existence condition which was found in \cite{KT} and concerned
only the topological equivalence of singular points. Another way
is a higher-dimensional generalization. For example, for
holomorphic function germs $f:(\C^n,0)\to(\C,0)$, it is known that
$$\sqrt[n]{\mu(f)}+1\le d(f)\le\mu(f)+1\ .$$
We conjecture that
$$d(f)\le a_n\sqrt[n]{\mu(f)}$$
with $a_n>0$ depending only on $n$, and we prove this for germs of
type $A_k$, $k\ge 1$ (Theorem \ref{t6}, section \ref{sec7}).

{\bf Idea of the proof}. Similarly to \cite{GLS,Sh2} we introduce
certain zero-dimensional schemes $Z\subset\PP^2$ associated with
singular and critical points, whose degree (length) is bounded by
a linear function of the (total) Milnor number, and such that, for
$$d=\min\{n\ge 1\ :\
H^1(\PP^2,\kj_Z(n))=0\},\quad\kj_Z=\Ker\left(
\ko_{\PP^2}\to\ko_Z\right)\ ,$$ there is
a curve (polynomial) of degree $d$ with singular (critical) points
of given types (Lemmas \ref{l5}(2), section \ref{secnew1}, Lemma
\ref{l8}, section \ref{secnew2},and proof of Theorem \ref{t5},
section \ref{sec3}). In principle, $d$ may be as large as $\deg
Z-1$. We, however, can choose $Z$ to be {\it generic} in
$\Iso(Z)$, the set of zero-dimensional schemes isomorphic to $Z$
as subschemes of $\PP^2$, and then establish our principal bound
(Proposition \ref{p3}, section \ref{sec5})
\begin{equation}d<\frac{4}{\sqrt{3}}\sqrt{\deg Z}-2\quad\text{as}\quad
\deg Z>2\ ,\label{enew1}\end{equation} which provides the main
estimates for the analytic order of a singular or critical point.

In \cite{GLS} an upper bound like (\ref{enew1}) is obtained for
irreducible zero-dimensional schemes of cluster type, generic in
their deformation class (which can be rather larger than the
isomorphism class). The proof was based on the so-called ``Horace
method" suggested by Hirschowitz \cite{H}. It consists in an
inductive procedure, where on each step one specializes a
zero-dimensional scheme (in its deformation class) on a given
line, then passes to the residue scheme. However, this approach
fails in our situation. The main obstacle (besides many technical
ones) is that, starting with a zero-dimensional scheme $Z$ generic
in $\Iso(Z)$, we have to specialize it in certain way, and then
obtain a residue scheme which is no longer generic in its
isomorphism class, thus, induction assumption does not apply.

To obtain (\ref{enew1}), we exploit a different idea, which is
similar in a sense to that in \cite{Xu1}, where $h^1$-vanishing
for some zero-dimensional schemes in the plane is deduced from the
ampleness of some divisors of the blown-up plane by Kodaira's
theorem. Namely, we start with estimating the minimal degree of a
curve, containing a scheme $Z$ generic in $\Iso(Z)$, from below by
$\sqrt{\deg Z}/2$ (Proposition \ref{p1}, section \ref{sec5}). For,
we observe that, deforming $Z$ in $\Iso(Z)$ so that the curve of
minimal degree through the scheme changes, one obtains that either
an intersection of $Z$ with a close element of $\Iso(Z)$ is of a
(relatively) large length, or $Z$ determines a singularity with
large invariants so that the curves through close schemes have
many intersections in neighborhood of singularities, and then the
desired estimate comes from B\'ezout's theorem. A combination of
such arguments can be found in \cite{Xu1} when zero-dimensional
schemes define ordinary singular points; in general case we use
estimates from \cite{GS}. Next, instead of exploring ampleness
which seems to be not easy to apply in our problem, we use the
Castelnuovo function theory (see \cite{D,GLS1}). The latter
argument appears to be quite simple and transparent. The graph of
the (positive) Castelnuovo function of a zero-dimensional scheme
$Z$ has width $d-1$, where $d$ is from (\ref{enew1}), its height
equals the minimal degree of a curve through $Z$, which is
$\sim\sqrt{\deg Z}$, and the area of its convex hull is $\deg Z$.
Thus, one obtains $d\sim\sqrt{\deg Z}$ when getting rid of long
horizontal segments of the graph. The latter can be done by Davis'
lemma \cite{D} (see details in the proof of Proposition \ref{p3},
section \ref{sec5}, and in \cite{GLS1}).

Finally, we notice that one could similarly treat zero-dimensional
schemes $Z\subset\PP^n$, $n\ge 3$. Indeed, it is not difficult to
show that the minimal degree of a hypersurface through a scheme
$Z$ generic in $\Iso(Z)$ is $\sim\sqrt[n]{\deg Z}$. However, the
lack of an appropriate Castelnuovo function theory prevents to
make step to $h^1$-vanishing bounds.

\begin{remark} One may ask what is the minimal possible
coefficient of $\sqrt{\deg Z}$ in (\ref{enew1}). An example of two
``fat"\footnote{That is defined by powers of the maximal ideals}
points of equal multiplicities shows that it cannot be less than
$2$ which is close to our value $4/\sqrt{3}=2.30...$

The Harbourne-Hirschowitz conjecture \cite{H} (see a survey and
bibliography in \cite{CM}) states that, for the scheme $Z$ of
``fat" points in general position in the plane,
$h^1(\PP^2,\kj_Z(d))=0$ as far as $h^0(\PP^2,\kj_Z(d))\ge 0$ and
$d$ is greater or equal to the sum of any three multiplicities,
i.e., $d\approx\sqrt{2\deg Z}$ in (\ref{enew1}). It is not clear
what should be an analogue of this conjecture for arbitrary
schemes if any. A reasonable conjecture can be an analogue of the
Alexander-Hirschowitz theorem \cite{AH}: for any $k\ge 1$ there
exists $N(k)\ge 1$ such that $h^1(\PP^2,\kj_Z(d))=0$, provided,
$d\ge N(k)$ and $h^0(\PP^2,\kj_Z(d))\ge 0$, for any scheme
$Z\subset\PP^2$ with irreducible components of length $\le k$ and
which is generic in $\Iso(Z)$.
\end{remark}

{\bf Acknowledgements.} I would like to thank G.-M. Greuel and C. Lossen for
very useful remarks and comments which allowed me to correct mistakes and
improve the presentation.

\section{Zero-dimensional schemes associated with
singular and critical points}

\subsection{Zero-dimensional schemes: cluster schemes, numerical
invariants, deformation and isomorphism classes} Throughout the
paper, we work with zero-dimensional schemes $Z$ that are
contained in a smooth algebraic surface $\Sig$. The  corresponding
ideal sheaves will be denoted by $\kj_{Z/\Sig}\subset\ko_{\Sig}$.
Moreover, we denote
$$\deg Z=\sum_z\dim_\C
\hat\ko_{\Sig,z}/(\kj_{Z/\Sig})_z\,, \qquad \mt(Z,z)=\max
\left\{n\in
  \Z\ :\ (\kj_{Z/\Sig})_z\subset \mm_z^n \right\}\,, $$
with $\hat\ko_{\Sig,z}$ the analytic local ring at $z$ and
$\mm_z\subset\hat\ko_{\Sig,z}$ the maximal ideal.

Let $z$ be an isolated singular point of an algebraic curve $C
\subset\Sig$. Denote by $T^{\infty}(C,z)$
the complete (infinite) embedded resolution tree of the singular point
$(C,z)$. The root vertex of $T^{\infty}(C,z)$ is $z$.
Among the other vertices of $T^{\infty}(C,z)$ there are finitely many
which are not nodes of the union of the exceptional locus with the
corresponding strict transform of $(C,z)$. All these vertices
of $T^{\infty}(C,z)$ together with $z$ are called
{\it essential infinitely near points} of $(C,z)$.
They are vertices of a subtree $T^*(C,z)\subset T^{\infty}(C,z)$.

Define the {\it multiplicity} of $C$ at $z$ as
$$\mt(C,z)=\max\{n\ge 0\ :\ C\in\mm_z^n\}\ .$$
Correspondingly, for any vertex $q\in T^{\infty}(C,z)$,
by the multiplicity $\mt(C,q)$ of $C$ at $q$ we mean the
multiplicity of the respective strict transform of $C$ at $q$.
For example,
$$\mt(C,q)=1,\quad q\in T^{\infty}(C,z)\backslash T^*(C,z)\ .$$

Given a finite subtree $T\subset T^{\infty}(C,z)$
such that $T\supset T^*(C,z)$, we define a {\it
zero-dimensional cluster scheme} $Z=Z(T)$ by the
ideal $I(Z)\subset\hat\ko_{\Sig,z}$ generated by
all germs $f\in\hat\ko_{\Sig,z}$ satisfying
$$\mt(f,q)=\mt(C,q)\stackrel{\text{\rm def}}{=}\mt(Z,q),\quad q\in T\ .$$
Such schemes were introduced in \cite{GLS} as
``generalized singularity schemes". The above generators
$f$ of the ideal $I(Z)$ are called generic elements of $I(Z)$.

\begin{lemma}\label{l3} {\rm (\cite{GLS}, Lemma 2.4)}
(1) The vertices of $T=T(Z)$ form the base point set of the ideal
$I(Z)$.

(2) Almost all germs $f\in I(Z)$ are generic.
\end{lemma}

A reducible zero-dimensional scheme $Z\subset\Sig$, concentrated
at points $z_1,...,z_p$, whose irreducible components
$Z_{z_1},...,Z_{z_p}$ are cluster schemes, is called a cluster
scheme as well, and $T(Z)$ is defined as the disjoint union of the
trees $T(Z_{z_i})$, $i=1,...,p$.

Let $Z$ be an arbitrary irreducible zero-dimensional scheme in
$\PP^2$. There exists a unique maximal cluster subscheme
$Z_{cl}\subset Z$. Namely, $T(Z_{cl})$ is the tree of infinitely
near base points of the ideal $I(Z)$, and $\mt(Z_{cl},q)$, $q\in
T(Z_{cl})$, are the corresponding multiplicities of a generic
element of $I(Z)$; more precisely, a generic element in the linear
space spanned by (finitely many) generators of $I(Z)$. If $Z$ is
reducible, then $Z_{cl}$ is the union of the maximal cluster
subschemes of the components of $Z$.

Put
$$M_2(Z)=\sum_{q\in T(Z_{cl})}(\mt(Z_{cl},q))^2\ .$$
We claim that
\begin{equation}
\deg Z\le M_2(Z)=\deg Z_{cl}+\del\ ,\label{e38}
\end{equation}
where $\del$ is the $\del$-invariant of a generic member of
$I(Z)$. The equality in (\ref{e38}) follows from the formulas in
\cite{GLS}, Lemma 2.6. For the inequality suppose that $Z$ is
concentrated at point $z$. Then take two distinct generic elements
$f,g\in I(Z)$, and obtain $$Z\subset f\cap
g\quad\Longrightarrow\quad \deg Z\le(f\cdot g)_z=\sum_{q\in
T(Z_{cl})}(\mt(Z_{cl},q))^2=M_2(Z) \ .$$ Note that (\ref{e38})
implies
\begin{equation}
M_2(Z)<2\cdot\deg Z\ .\label{e73}
\end{equation}

For a zero-dimensional scheme $Z$ with $p$ irreducible components,
its isomorphism class $\Iso(Z)$ is fibred over the space of
$p$-tuples $(z_1,...,z_p)\in(\PP^2)^p$ with fibre being an orbit
of the action of the group
$\prod_{i=1}^p\Aut(\hat\ko_{\PP^2,z_i}/\mm_{z_i}^n)$, where $n$ is
sufficiently large. If $Z$ is a cluster scheme then one can
naturally define the set $\Def(Z)$ of schemes $Z'\subset\PP^2$,
deformation equivalent to $Z$. This is a smooth irreducible
quasiprojective variety \cite{GLS1}, section 2. Clearly,
$\Def(Z)\supset\Iso(Z)$.

A zero-dimensional scheme is called {\it nonsingular}, if the
generic elements of the ideals of its components are nonsingular,
and is called {\it singular} otherwise. We notice that a
nonsingular zero-dimensional scheme is always a cluster scheme.

\subsection{Zero-dimensional schemes associated with topological types of
singular points}\label{secnew1} Let $z$ be an isolated singular
point of a curve $C\subset\Sig$. The {\it equisingularity ideal}
introduced in \cite{Wa} (see also \cite{DH,GL,GLS,GLS1,Sh}) is
defined as
$$I^{es}(C,z)\,:=\,\bigl\{\, g\in \hat\ko_{\Sig,z}\ :\ f+\varepsilon g \text{
  is equisingular over Spec}\,(\C[\varepsilon]/\varepsilon^2)\, \bigr\}\ ,$$
where $f\in\hat\ko_{\Sig,z}$ is a germ induced by $C$. The
zero-dimensional scheme
defined by $I^{es}(C,z)$ is denoted by $Z^{es}(C,z)$.

The notion of the cluster scheme directly relates to the
topological equivalence of germs. Namely, \cite{GLS}, Lemma 2.4,
implies

\begin{lemma}\label{l4}
For any irreducible deformation equivalent cluster schemes $Z$ and $Z'$,
generic elements $f\in I(Z)$, $f'\in I(Z')$ are topologically equivalent.
\end{lemma}

Let $(C,z)\subset(\PP^2,z)$ be a reduced curve germ. We associate three
zero-dimensional schemes with it (cf. \cite{GLS}):
\begin{itemize}
\item $Z^s(C,z)$, the cluster scheme defined by the germ $(C,z)$ and
the tree of essential point $T^*(C,z)$;
\item $Z^s_1(C,z)=Z^s(\widetilde C,z)$, the cluster scheme defined
by the germ $(LC,z)$, where $L$ is a straight line through $z$ and
transverse to $C$.
\end{itemize}

The importance of these schemes arises from

\begin{lemma}\label{l5}
(1) The scheme $Z^s(C,z)$ is minimal among the zero-dimensional schemes
$Z$ such that almost all germs $f\in I(Z)$ are topologically equivalent to
$(C,z)$.

(2) Let some zero-dimensional scheme $Z$ such that $z\not\in Z$,
satisfy
\begin{equation}
H^1(\kj_{Z^s_1(C,z)\cup Z/\PP^2}(d))=0\ .\label{e37}
\end{equation}
Then there exist $D\in |\kj_{Z^s(C,z)\cup Z/\PP^2}(d)|$ such that
$(D,z)$ is topologically equivalent to $(C,z)$. Moreover, these
curves $D$ form a dense open subset in $|\kj_{Z^s(C,z)\cup
Z/\PP^2}(d)|$.
\end{lemma}

Clearly $Z^s(C,z)\supset Z^{es}(C,z)$. Moreover, it can be shown that
$Z^s(C,z)$ is the minimal cluster scheme containing $Z^{es}(C,z)$.

{\it Proof}.
The first statement reflects the fact that the
tree of essential infinitely
near points of $(C,z)$ and multiplicities of $C$ at them
(uniquely) determine the topological type of $(C,z)$.

The second statement can be proven in the same way as it is done
in Step 1 of the proof of Lemma 5.8 in \cite{GLS}, where $Z$ is
supposed to be empty, $Z_1^S(C,z)$ is denoted by $\widetilde X$,
and the required $h^1$-vanishing condition is found in (5.12).
Note only that the scheme $X'$, used in this proof, is a subscheme
of $\widetilde X$; hence the $h^1$-vanishing for $X'$ mentioned in
(5.12) follows from that for $\widetilde X=Z_1^s(C,z)$. \proofend

\subsection{Zero-dimensional schemes associated with analytic types of
singular points}\label{secnew2} In the previous notation, we
introduce the zero-dimensional scheme $Z^{ea}(C,z)$ defined by the
{\it Tjurina ideal}
$$I^{ea}(C,z)\,:=\, \langle\, f, f_x,
f_y\,\rangle \,\subset\,
\hat\ko_{\Sig,z}\,,$$ where $f(x,y)=0$ is a local equation for the
germ $(C,z)$. The ideal $I^{ea}(C,z)$ is the tangent space to
equianalytic (i.e., analytically trivial) deformations of $(C,z)$.

For our purpose we shall use zero-dimensional schemes associated
with analytic types of singular points, which are analogous to
$Z^s,Z^s_1$, but without the minimality property as in Lemma
\ref{l5}(1) (except for the simple singularities $A_k$, $k\ge 1$,
$D_k$, $k\ge 4$, $E_6$, $E_7$, $E_8$, for which topological and
analytic equivalence coincide).

If the singular point $z$ of
$C$ is simple, we put
$$Z^a(C,z)=Z^s(C,z),\quad Z^a_1(C,z)=Z^s_1(C,z)\ .$$
If the singular point $z$ of $C$ is not simple, $f(x,y)=0$ is an
equation of $C$ in a neighborhood of $z$, following \cite{GLS1},
section 1.3, and we define
\begin{itemize}
\item the zero-dimensional scheme $Z^a(C,z)$ by the ideal
$$I^a(C,z)=\left\{
g\in\hat\ko_{\PP^2,z}\ :\ \langle g,g_x,
g_y\rangle\subset
\langle f,f_x,
f_y\rangle\right\}\ ,$$
\item the zero-dimensional scheme $Z^a_1(C,z)$ by the
ideal $I^a_1(C,z)=\mm_z\cdot I^a(C,z)$.
\end{itemize}

Observe that $Z^a_1(C,z)\supset Z^a(C,z)\supset Z^{ea}(C,z)\supset Z^{es}(C,z)$.

\begin{lemma}\label{l7}
In the above notation,

(1) any scheme $Z\in\Iso(Z^a(C,z))$
(resp., $Z\in\Iso(Z^a_1(C,z))$) is $Z^a(\widetilde C,w)$
(resp., $Z^a_1(\widetilde C,w)$)
for some germ $(\widetilde C,w)$ analytically equivalent to
$(C,z)$;

(2) if $g\in I^a(C,z)$ (resp., $g\in I^a_1(C,z)$), then, for
almost all $t\in\C$ (resp., for all $t\in\C$) the curve germs
$\{f=0\}$ and $\{f+tg=0\}$ are analytically equivalent.
\end{lemma}

{\it Proof}. The first statement is evident. The fact that, for
$g\in I^a(C,z)$ and almost all $t\in\C$, the germs $\{f=0\}$ and
$\{f+tg=0\}$ are analytically equivalent, follows from the
Mather-Yau theorem \cite{MY} (see also Lemma 1.8(a) \cite{GLS1}).

Assume now that $g\in I^a_1(C,z)=\mm_z\cdot I^a(C,z)$. Observe
that if $h\in I^a(C,z)$, then $h=af+bf_x+cf_y$, $b,c\in\mm_z$,
which follows from Lemma 1.8(a,c) \cite{GLS1}. Then
\begin{eqnarray}
&g=af+(xb'+yb'')f_x+(xc'+yc'')f_y,\quad
a,b',b'',c',c''\in\mm_z,\label{e400}\\
&b'f_{xx}+c'f_{xy},\ b''f_{xx}+c''f_{xy},\ b'f_{xy}+c'f_{yy}, \
b'f_{xy}+c'f_{yy}\in\langle f,f_x,f_y\rangle\ .\label{e401}
\end{eqnarray}
Since $1+a\in\hat\ko_{\PP^2,z}^*$, the germ
$f+g=f+af+(xb'+yb'')f_x+(xc'+yc'')f_y$ is equivalent to
$f+\widetilde g$,
$$\widetilde g=\frac{xb'+yb''}{1+a}f_x+\frac{xc'+yc''}{1+a}f_y\
.$$ The restrictions to $b',b'',c',c''$ in (\ref{e400}),
(\ref{e401}) yield that
$$\langle\widetilde g,\widetilde g_x,\widetilde g_y\rangle\subset
\mm_z\cdot\langle f,f_x,f_y\rangle\ ,$$ so the equivalence of $f$
and $f+\widetilde g$ follows from the Mather-Yau theorem \cite{MY}
(cf. Lemma 1.8(b) \cite{GLS1}).

\begin{lemma}\label{l8}
Let a zero-dimensional scheme $Z$ such that $z\not\in Z$ satisfy
\begin{equation}
H^1(\kj_{Z_1^a(C,z)\cup Z/\PP^2}(d))=0\ .\label{e49}
\end{equation}
Then there exist a curve $D\in |\kj_{Z^a(C,z)\cup Z/\PP^2}(d)|$,
whose germ $(D,z)$ is analytically equivalent to $(C,z)$.
Moreover, such curves $D$ form a dense open subset in
$|\kj_{Z^a(C,z)\cup Z/\PP^2}(d)|$.
\end{lemma}

{\it Proof}. In the exact sequence
$$H^0(\kj_{Z/\PP^2}(d))\to \hat\ko_{\PP^2,z}/I^a_1(C,z)
\to H^1(\kj_{Z\cup Z^a_1(C,z)/\PP^2}(d))=0\ ,$$ the first morphism
is surjective. Denote by $\varphi\in\hat\ko_{\PP^2,z}/I^a_1(C,z)$
the image of a germ $\psi\in\hat\ko_{\PP^2,z}$ defined by $(C,z)$.
Take $\Phi\in H^0(\kj_{Z/\PP^2}(d))$ which projects to $\varphi$.
Then $\Phi-\psi\in I^a_1(C,z)$, and by Lemma \ref{l7}(2), the
curve germs $(\{\Phi=0\},z)$ and $(C,z)$ are analytically
equivalent, thus, we can put $D=\{\Phi=0\}$.

\subsection{Zero-dimensional schemes associated with analytic types
of critical points}\label{sec2} Let $f:(\C^2,0)\to(\C,0)$ be a
germ of a holomorphic function with a finite Milnor number
$\mu(f)=\dim\hat\ko_{\C^2,0}/\langle f_x,f_y\rangle$. Germs
$f,g\in\hat\ko_{\C^2,0}$ are called {\it (right) equivalent}, if
there is $\psi\in\Aut(\hat\ko_{\C^2,0})$ such that $g=f\circ\psi$.
Introduce the zero-dimensional schemes
 \begin{itemize}
 \item $Z_0(f)$ defined by the ideal
 $$I_0(f)=\{g\in\hat\ko_{\C^2,0}\ :\ g,g_x,g_y\in
\langle f_x,f_y\rangle\}\ ;$$
 \item $Z(f)$ defined by the
ideal $I(f)=\mm_0I_0(f)$.
 \end{itemize}
 An analogue of Lemma \ref{l7} reads as

\begin{lemma}\label{l11}
In the above notation,

(1) any scheme $Z\in\Iso(Z(f))$
is $Z(\widetilde f)$
for some germ $\widetilde f$ equivalent to $f$;

(2) if $g\in I(f)$, then, for all $t\in\C$, the germ $f+tg$ is
equivalent to $f$.
\end{lemma}

The first statement of Lemma is evident. The second one is, in
fact, known and can be proven as Mather's finite determinacy
theorem \cite{Ma}.

\subsection{Bounds for degrees of zero-dimensional schemes}
Given a reduced curve germ $(C,z)\subset\PP^2$ of a function germ
$f\in\hat\ko_{\C^2,0}$, the degrees of the schemes $Z^s(C,z)$,
$Z^a(C,z)$, $Z(f)$ are invariants of the given singular or
critical point up to the corresponding equivalence. We shall
compare these invariants with the classical ones.

\begin{lemma}\label{l10}
In the above notation,
\begin{equation}
\deg Z^s(C,z)=\deg Z^a(C,z)=\deg
Z_0(f)=\left[\frac{3k+4}{2}\right]\ ,\label{e41}
\end{equation}
if $(C,z)$ or $f\in\hat\ko_{\C^2,0}$ is of type $A_k$, $k\ge 1$.
For other singular and critical points
\begin{eqnarray}
&\deg Z^s(C,z)\le 3\del(C,z)\ ,\label{e40}\\ &\deg Z^a(C,z)\le
2\mu(C,z)\ ,\label{e42}\\ &\deg Z_0(f)\le 3\mu(f)-2\cdot\mt(f)+2\
.\label{e71}
\end{eqnarray}
\end{lemma}

{\it Proof}. Formula (\ref{e41}) is an easy computation along the
definition of $Z^s(C,z)=Z^a(C,z)=Z_0(f)$ in this case.

Assume that $\mt(C,z)\ge 3$.

Then
$$\deg Z^s(C,z)=\sum_{q\in T^*(C,z)}\frac{\mt(C,q)\cdot(\mt(C,q)+1)}{2}
=\del(C,z)+\sum_{q\in T^*(C,z)}\mt(C,q)\ .$$
By definition of $T^*(C,z)$,
$$\#\{q\in T^*(C,z)\ :\ \mt(C,q)=1\}\le\mt(C,z)=m\ .$$
Hence
$$\sum_{q\in T^*(C,z)}\mt(C,q)\le 2\cdot\mt(C,z)+
\sum_{\renewcommand{\arraystretch}{0.6}
\begin{array}{c}
\scriptstyle{q\in T^*(C,z)}\\
\scriptstyle{q\ne z}\\
\scriptstyle{\mt(C,q)>1}
\end{array}}\mt(C,q)$$
$$\le 2\sum_{q\in T^*(C,z)}\frac{\mt(C,z)\cdot
(\mt(C,z)-1)}{2}=2\del(C,z)\ ,$$
and (\ref{e40}) follows.

We shall establish (\ref{e42}), first, for simple singularities
$D_k$, $k\ge 4$, $E_k$, $k=6,7,8$. Here $Z^a(C,z)=
Z^s(C,z)$, and a direct computation gives
$$\deg Z^a(D_k)=\deg Z^s(D_k)=\left[\frac{3k+1}{2}\right]\le 2k,
\quad k\ge 4\ ,$$
$$\deg Z^a(E_k)=\deg Z^s(E_k)=k+3\le 2k,\quad k=6,7,8\ .$$
If $(C,z)$ is not simple, introduce $\Pi_1,\Pi_2$,
two distinct generic polar curves of $C$,
and $\Pi_{11},\Pi_{12}$,
two generic polar curves of $\Pi_1$. By \cite{GLS1}, formula (1.5),
$$I^a(C,z)\supset\{fC+g\Pi_1\ :\ f,g
\in\hat\ko_{\PP^2,z},\ g\Pi_{11},g\Pi_{12}\in\langle\Pi_1,\Pi_2\rangle\}\ .$$
By the double point divisor theorem \cite{VDW}, \S 50, the ideal
$$\{g\in\hat\ko_{\PP^2,z}\ :\ g\Pi_{11},g\Pi_{12}\in\langle\Pi_1,
\Pi_2\rangle\}$$
contains the ideal
$$I=\{g\in\hat\ko_{\PP^2,z}\ :\ (g\cdot P)_z\ge(\Pi_2\cdot P)_z-\mt(P,z)+1$$
$$\text{for any local branch}\ P\ \text{of}\ (\Pi_1,z)\}\ ,$$
where $(*,*)_z$ denotes the intersection multiplicity of two curve germs
at the point $z$\footnote{For the sake of notation we write $g$ instead of
$\{g=0\}$ in these formulas.}.
Here
$$\dim\hat\ko_{\PP^2,z}/I=\sum_P((\Pi_2\cdot P)_z-\mt(P,z)+1)-\del(\Pi_1,z)$$
$$\le\sum_P(\Pi_2\cdot P)_z-\del(\Pi_1,z)=\mu(C,z)-\del(\Pi_1,z)\ .$$
If $\mt(C,z)=3$ and $(C,z)$ is not simple, then $\Pi_1$ has at
least a tacnode at $z$, so $\del(\Pi_1,z)\ge 2=\mt(C,z)-1$. If
$\mt(C,z)\ge 4$, then $\mt(\Pi_1,z)= \mt(C,z)-1\ge 3$, so
$$\del(\Pi_1,z)\ge\frac{\mt(\Pi_1,z)\cdot(\mt(\Pi_1,z)-1)}{2}\ge\mt(\Pi_1,z)
=\mt(C,z)-1\ ,$$
which altogether results in $\dim\hat\ko_{\PP^2,z}/I\le
\mu(C,z)-\mt(C,z)+1$.
Hence
$$\deg Z^a(C,z)=\dim\hat\ko_{\PP^2,z}/I^a(C,z)\le
(C\cdot\Pi_1)_z+\dim\hat\ko_{\PP^2,z}/I$$
$$\le(\mu(C,z)+\mt(C,z)-1)+(\mu(C,z)-\mt(C,z)+1)=2\mu(C,z)\ .
$$

For inequality (\ref{e71}) we note that
 $$I_0(f)\supset\{g\in\hat\ko_{\C^2,0}\ :\ g=af_x+bf_y,\
 a,b\in\hat\ko_{\C^2,0}$$
 $$\qquad af_{xx},af_{xy},bf_{xy},bf_{yy}\in\langle
 f_x,f_y\rangle\}$$
 $$\supset\{g\in\hat\ko_{\C^2,0}\ :\ g=af_x+bf_y,\
 a,b\in I\}\ ,$$
 where the ideal $I\subset\hat\ko_{\C^2,0}$ is defined as in the
 preceding paragraph. Hence as in the previous computation
 $$\deg Z_0(f)\le\dim\hat\ko_{\C^2,0}/\langle f_x,f_y\rangle+
 2\dim\hat\ko_{\C^2,0}/I\le 3\mu(f)-2\cdot\mt(f)+2\ .$$

\section{Analytic and
topological order of a zero-dimensional scheme
in the plane}

\subsection{Definitions and notations}\label{sec5}
A zero-dimensional scheme $Z\subset\PP^2$ can be characterized by
the following numbers (orders)
\begin{eqnarray}
&\ord_0(Z)=\min\{d\in\Z\ :\ H^0(\kj_{Z/\PP^2}(d))>0\}\ ,\nonumber\\
&\ord_1(Z)=\min\{d\in\Z\ :\ H^1(\kj_{Z/\PP^2}(d))=0\}\ .\nonumber
\end{eqnarray}
Put
$$\ord^{top}_0(Z)=\max_{Z'\in\Def(Z)}\ord_0(Z'),\ \ord^{top}_1(Z)
=\min_{Z'\in\Def(Z)}\ord_1(Z'),\ Z\ \mbox{is
cluster scheme},$$
$$\ord^{an}_0(Z)=\max_{Z'\in\Iso(Z)}\ord_0(Z'),\quad
\ord^{an}_1(Z)=\min_{Z'\in\Iso(Z)}\ord_i(Z'),\quad Z\ \mbox{is any
scheme}.$$ Clearly, $\ord^{top}_i(Z)$, $\ord^{an}_i(Z)$ are
$\ord_i(Z')$ for a generic element $Z'$ of the corresponding
family.

It was shown in \cite{GLS}, Lemmas 3.1, 4.1, 5.8, that
\begin{itemize}
\item for an irreducible zero-dimensional cluster scheme $Z$
defined by a curve germ having only nonsingular local
branches,
\begin{equation}
\ord^{top}_1(Z)<(1+\sqrt{2})\sqrt{\deg Z}+\mt Z+1\ .\label{e19}
\end{equation}
\item for an arbitrary irreducible zero-dimensional cluster scheme $Z$,
\begin{equation}
\ord^{top}_1(Z)<\frac{\sqrt{85}-3}{2}\sqrt{\deg Z}+\mt Z+
\mt_sZ+1\ ,\label{e20}
\end{equation}
where $\mt_sZ$ is the sum of multiplicities of the singular branches
of the germ $f(Z)$.
\end{itemize}

We shall estimate the analytic orders of any zero-dimensional
scheme. As a by-product we improve estimates (\ref{e19}),
(\ref{e20}) for topological orders and extend them to reducible
schemes.

\begin{proposition}\label{p1}
For any zero-dimensional scheme $Z\subset\PP^2$,
\begin{equation}
\ord^{an}_0(Z)\ge\frac{\deg Z}{\sqrt{2M_2(Z)}} \ .\label{e3}
\end{equation}
\end{proposition}

\begin{remark}\label{r3}
In view of (\ref{e73}), the bound (\ref{e3}) can be weakened up to
the following, simpler inequality:
$$\ord^{an}_0(Z)>\frac{\sqrt{\deg Z}}{2}\ .$$
\end{remark}

{\it Proof}. {\bf Step 1}. Consider, first, the case of
an irreducible scheme $Z$ concentrated at a point
$z\in\PP^2$.

Let $Y$ be a generic element of $\Iso(Z)$,
concentrated at $z$,
and $d=\ord_0(Y)$. Take a generic curve
$C\in|\kj_{Y/\PP^2}(d)|$, and suppose that
$C=C_1^{l_1}...C_r^{l_r}$, where $C_1,...,C_r$ are
distinct reduced irreducible curves of degrees
$d_1,...,d_r$, respectively, so that
$d=l_1d_1+...+l_rd_r$.

If $C$ contains only point $z$ of $T(Y_{cl})$, then $C$
transversally intersects a generic element $f\in I(Y)$. Since
$d\ge\mt(Y_{cl},z)$, we have
$$d\ge\sqrt{d\cdot\mt(Y_{cl},z)}=\sqrt{(C\cdot f)_z}
\ge\sqrt{\deg Y}=\sqrt{\deg Z}\ .$$

Now without loss of generality we can assume that the
tree $T(Y_{cl})$ contains at least two points, the curves
$C_i$, $1\le i\le s<r$, contain only point $z$ of $T(Y_{cl})$,
and any curve $C_i$, $i>s$, contains at least two
points of $T(Y_{cl})$.

Fix $s<i\le r$. Put $T_i=T(C_i)\cap T(Y_{cl})$. Choose local
coordinates $x,y$ in a neighborhood of $z=(0,0)$
so that the axes are transverse to $C_i$ and to a generic $f\in I(Y)$.
Introduce the sequence
\begin{equation}
\psi_m\in\Aut(\hat\ko_{\PP^2,z}),\quad\psi_m(x,y)=(x,y+\eps x^m),
\quad m\ge 1,\quad\eps=\const\ne 0\ ,\label{e11}
\end{equation}
and consider
the schemes $\psi_m(Y)$. The ascending sequence of trees
$T(Y_{cl})\cap T(\psi_m(Y_{cl}))$ stabilizes for a sufficiently large $m$.
Then there exists
\begin{equation}
k=\max\{m\ge 1\ :\ T_i\not\subset T(\psi_m(Y_{cl}))\}\ .\label{e26}
\end{equation}

\begin{lemma}\label{l1}
\begin{equation}
\sum_{q\in T_i\backslash T(\psi_k(Y_{cl}))}
(\mt(C_i,q))^2\le(\mt(C_i,z))^2\ .\label{e4}
\end{equation}
\end{lemma}

{\it Proof of Lemma \ref{l1}}.
Denote by $T_m(C_i)$, $m\ge 1$, the tree of common
infinitely near points of the curves $C_i$ and
$\psi_m(C_i)$ at $z$. Then
$$T_i\backslash T(\psi_k(Y_{cl}))\subset T_{k+1}(C_i)
\backslash T_k(C_i)$$
$$\Longrightarrow\quad
\sum_{q\in T_i\backslash T(\psi_k(Y_{cl}))}
(\mt(C_i,q))^2\le
\sum_{q\in T_{k+1}(C_i)\backslash T_k(C_i)}(\mt(C_i,q))^2$$
$$=(C_i\cdot\psi_{k+1}C_i)_z-(C_i\cdot\psi_kC_i)_z\ .$$
To estimate the latter expression, we use the Puiseux
decomposition
\begin{equation}
C_i(x,y)=(1+O(x,y))\prod_{s=1}^n(y-\xi_s(x)),\quad n=\mt(C_i,z),\label{e28}
\end{equation}
where $\xi_s(x)$, $s=1,...,n$, are fractional power series.
Then
$$(C_i\cdot\psi_k(C_i))_z=\sum_{1\le r,s\le n}\ord(\xi_r(x)-\xi_s(x)-
\eps x^k)\ ,$$
$$(C_i\cdot\psi_{k+1}(C_i))_z=\sum_{1\le r,s\le n}\ord(\xi_r(x)-\xi_s(x)-
\eps x^{k+1})\ ,$$
and (\ref{e4}) follows, because of an obvious inequality
$$\ord(\xi_r(x)-\xi_s(x)-\eps x^{k+1})\le
\ord(\xi_r(x)-\xi_s(x)-\eps x^k)+1\ ,$$ where $\ord(*)$ means the
minimal power of $x$ occurring in the series $*$. \proofend

In view of the generality of $Y\in\Iso(Z)$ and $\eps\ne 0$,
the scheme $Y'=\psi_k(Y)$ is also generic in $\Iso_z(Z)$.
Hence $\ord_0(Y')=\ord_0(Y)=d$, and $Y'$ is contained in a curve
$C'$ of degree $d$ which splits as
$C'=(C'_1)^{l_1}...(C'_r)^{l_r}$, where $\deg C'_i=\deg C_i=d_i$,
and the curve $C'_i$ is different from $C_i$,
since $\psi_k$ moves the points
$q\in T_i\backslash T(Y'_{cl})\ne\emptyset$.
Then
$$d_i^2\ge\sum_{q\in T(Y_{cl})\cap T(Y'_{cl})}(\mt(C_i,q))^2\ ,$$
which by (\ref{e4}) implies
\begin{equation}
d_i^2\ge\frac{1}{2}\sum_{q\in T(Y_{cl})}(\mt(C_i,q))^2\ .
\label{e2}
\end{equation}
Note that (\ref{e2}) holds for $i=1,...,s$ as well.

Take a generic element $f\in I(Y)$.
Then $Y$ is contained in the scheme-theoretic intersection
of $f$ and $C$ at the point $z$. Hence
$$\deg Y=\deg Z\le(f\cdot C)_z=
\sum_{q\in T(Y_{cl})}\mt(Y_{cl},q)\cdot\mt(C,q)$$
$$\le\sqrt{\sum_{q\in T(Y_{cl})}(\mt(Y_{cl},q))^2}
\sqrt{\sum_{q\in T(Y_{cl})}(\mt(C,q))^2}
=\sqrt{M_2(Z)}\sqrt{\sum_{q\in T(Y_{cl})}(\mt(C,q))^2}\ .$$
Here
$$\sum_{q\in T(Y_{cl})}(\mt(C,q))^2=
\sum_{q\in T(Y_{cl})}\left(\sum_{i=1}^rl_i\cdot\mt(C_i,q)\right)^2
\stackrel{\text{(\ref{e2})}}{\le}2\left(\sum_{i=1}^rl_id_i\right)^2
=2d^2\ ,$$
and (\ref{e3}) follows.

{\bf Step 2}. Let $Z$ consist of components $Z_1,...,Z_p$
concentrated at points $z_1,...,z_p\in\PP^2$, respectively,
$p\ge 2$. Consider a generic scheme $Y\in\Iso(Z)$,
concentrated at a generic $p$-tuple $w_1,...,w_p\in\PP^2$,
and a generic curve $C\in|\kj_{Y/\PP^2}(d)|$,
$d=\ord_0(Y)=\ord^{an}_0(Z)$. Assume that
$C=C_1^{l_1}...C_r^{l_r}$, where $C_1,...,C_r$ are
reduced irreducible.

If a component $C_i$, $1\le i\le r$, passes through only one of
the points $w_1,...,w_p$, then (\ref{e2}) holds due to the argument in
Step 1 of the proof.

Let $C_i$ pass through points $w_1,...,w_s$, $s\ge 2$, and
be transverse to generic elements $f_j\in I(Y_{w_j})$,
$j=1,...,s$. Assuming that
$0<\mt(C_i,w_1)\le...\le\mt(C_i,w_k)$, we move the point $w_1$
to $w'_1\not\in C_i$ keeping $w_2,...,w_p$ fixed. The correspondingly
deformed scheme $Y'$ is also generic in $\Iso(Z)$, and hence there exists
a curve $C'\in|\kj_{Y'/\PP^2}(d)|$ splitting as
$C'=(C'_1)^{l_1}...(C'_r)^{l_r}$ with $\deg C'_j=d_j$, $j=1,...,r$,
$C'_i$ close to $C_i$ having multiplicities
$$\mt(C'_i,w'_1)=\mt(C_i,w_1),\quad\mt(C'_i,w_j)=\mt(C_i,
w_j),\ j=2,...,k\ .$$
Then by \cite{GS}, Theorem 2(1), or \cite{Xu}, Lemma 3,
$$d_i^2\ge(\mt(C_i,w_1))^2-\mt(C_i,w_1)+
\sum_{j=2}^k(\mt(C_i,w_j))^2$$
\begin{equation}
\ge\frac{1}{2}\sum_{j=1}^k(\mt(C_i,w_j))^2=\frac{1}{2}
\sum_{q\in T(Y_{cl})}(\mt(C_i,q))^2\ .\label{e10}
\end{equation}

Let $C_i$ pass through $w_1,...,w_s$, $s\ge 2$, and contain at
least two points of $T(Y_{w_1,cl})$. Then we apply transformations
(\ref{e11}) to $Y_{w_1}$, where $w_1=(0,0)$, keeping $Y_{w_j}$,
$j=2,...,p$, unchanged. The reasoning, similar to that in Step 1
of the proof, shows that a suitable $\psi_m$ moves the tree
$T(Y_{w_1,cl}) \cap T(C_i)$, turning $C_i$ into a curve $C'_i$ of
the same degree $d_i$ with
$$(C_i\cdot C'_i)_{w_1}\ge\frac{1}{2}\sum_{q\in T(Y_{w_1,cl})\cap C_i}
(\mt(C_i,q))^2\ ,$$
which immediately implies
$$d^2_i\ge\frac{1}{2}\sum_{q\in T(Y_{w_1,cl})\cap C_i}
(\mt(C_i,q))^2+\sum_{j=2}^k\sum_{q\in T(Y_{w_j,cl})\cap C_i}
(\mt(C_i,q))^2$$
\begin{equation}
\ge\frac{1}{2}\sum_{q\in T(Y_{cl})\cap C_i}(\mt(C_i,q))^2\ .\label{e12}
\end{equation}

Finally, (\ref{e3}) follows from (\ref{e2}), (\ref{e10}),
(\ref{e12}), as was done in Step 1 of the proof. \proofend

\begin{proposition}\label{p3}
For any singular zero-dimensional scheme $Z\subset\PP^2$,
\begin{equation}
\ord^{an}_1(Z) \le\sqrt{\frac{3}{2}M_2(Z)}+\frac{\deg
Z}{\sqrt{3M_2(Z)/2}}-2\ .\label{e7}
\end{equation}
For any nonsingular zero-dimensional cluster scheme $Z\subset\PP^2$,
\begin{equation}
\ord^{an}_1(Z)=-\left[\frac{3-\sqrt{1+8\deg Z}}{2}\right]\ .
\label{e22}
\end{equation}
\end{proposition}

\begin{remark}\label{r4}
Due to (\ref{e73}), inequality (\ref{e7}) implies a weaker
relation
\begin{equation}
\ord^{an}_1(Z) <\frac{4}{\sqrt{3}}\sqrt{\deg Z}-2\ ,\label{e24}
\end{equation}
which in turn is stronger than both (\ref{e19}) and (\ref{e20}).
\end{remark}

{\it Proof}. {\bf Step 1}. Let $Z$ be nonsingular (and hence a
cluster scheme).
Relation (\ref{e22}) means that
$$d=\ord^{an}_1(Z)=\min\{n\ :\ h^0(\ko_{\PP^2}(n))=
\frac{(n+1)(n+2)}{2}\ge\deg Z\}\ .$$ In other words, the $\deg Z$
conditions imposed on curves $C\in|\kj_{Y/\PP^2}(d)|$ by a generic
element $Y\in\Def(Z)=\Iso(Z)$ are independent. Take a natural
sequence of schemes $\emptyset\subsetneq
Z_1\subsetneq...\subsetneq Z_s=Z$, $s=\deg Z$, where the linear
system $|\kj_{Z_{i+1}/\PP^2}(d)|$ is obtained from
$|\kj_{Z_i/\PP^2}(d)|$ by imposing one condition, either a passage
through a point $\not\in Z_i$, or an extended by $1$ tangency
order with a fixed nonsingular curve germ. Then one can
inductively show that generic members of these linear systems are
nonsingular, and a generic choice of the new condition reduces the
dimension each time by $1$.

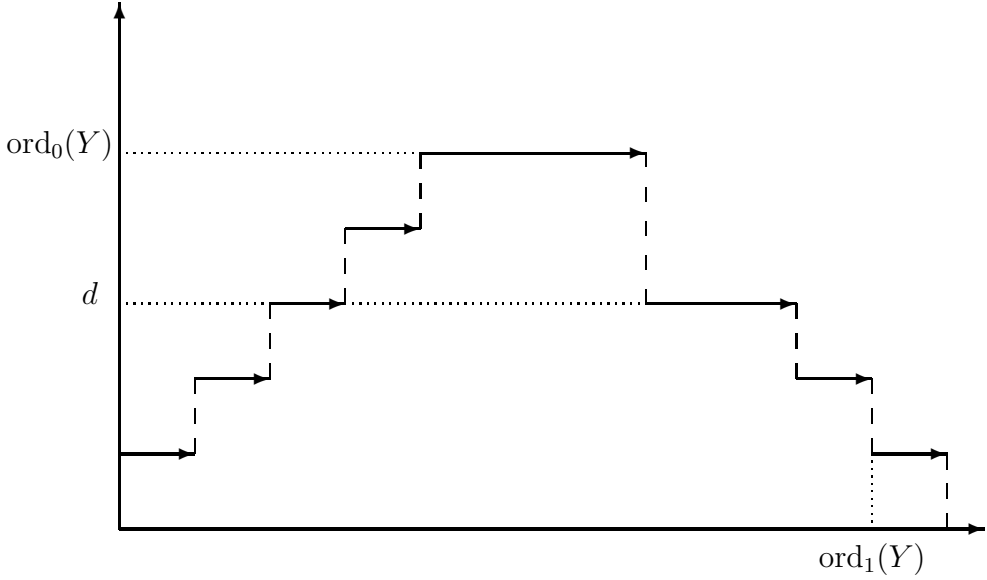
\begin{figure}
\setlength{\unitlength}{1cm}
\begin{picture}(13,8)(0,0)
\thicklines \put(2,1){\vector(1,0){11.5}}
\put(2,1){\vector(0,1){7}} \put(2,2){\vector(1,0){1}}
\put(3,3){\vector(1,0){1}} \put(4,4){\vector(1,0){1}}
\put(5,5){\vector(1,0){1}} \put(6,6){\vector(1,0){3}}
\put(9,4){\vector(1,0){2}} \put(11,3){\vector(1,0){1}}
\put(12,2){\vector(1,0){1}} \thinlines \dashline{0.2}(3,2)(3,3)
\dashline{0.2}(4,3)(4,4) \dashline{0.2}(5,4)(5,5)
\dashline{0.2}(6,5)(6,6) \dashline{0.2}(9,6)(9,4)
\dashline{0.2}(11,4)(11,3) \dashline{0.2}(12,3)(12,2)
\dashline{0.2}(13,2)(13,1) \dottedline{0.1}(2,6)(6,6)
\dottedline{0.1}(2,4)(4,4) \dottedline{0.1}(5,4)(9,4)
\dottedline{0.1}(12,2)(12,1) \put(0.5,6){$\ord_0(Y)$}
\put(1.5,4){$d$} \put(11.3,0.5){$\ord_1(Y)$}
\end{picture}
\caption{Castelnuovo function}\label{f1}
\end{figure}

{\bf Step 2}. Take a generic scheme $Y\in\Iso(Z)$. Consider the
Castelnuovo function
$\kc_Y(n)=h^1(\kj_{Y/\PP^2}(n-1))-h^1(\kj_{Y/\PP^2}(n))$, $n\ge
0$. One can find a detailed description of the Castelnuovo
function and its graph in \cite{D, GLS1}. In particular (see
Figure \ref{f1}),
\begin{eqnarray}
&\kc_Y(n)=n+1,\ 0\le n<\ord_0(Y)=\ord^{an}_0(Z)\ ,\label{e16}\\
&\kc_Y(n)\le\kc_Y(n-1),\ n\ge\ord_0(Y)\ ,\label{e13}\\
&\kc_Y(n)=0,\ n>\ord_1(Y)=\ord^{an}_1(Z)\ ,\label{e17}\\
&\sum_{n\ge 0}\kc_Y(n)=\deg Z=\deg Y\ .\label{e14}
\end{eqnarray}

If there is no $n\ge\ord_0(Y)$ such that $\ord_0(Y)>
\kc_Y(n)=\kc_Y(n-1)>0$, then we put $d=\ord_0(Y)$. If there exists
$0<d<\ord_0(Y)$ such that $\kc_Y(n)=\kc_Y(n-1)=d$ for some
$n\ge\ord_0(Y)$ (so called ``long stair"), then we assume that $d$
is minimal with this property (see Figure \ref{f1}). By \cite{D},
Claims 2.2 and 2.3, (see also \cite{GLS1}), there exists a curve
$C$ of degree $d$ such that (see Figure \ref{f1})
\begin{equation}
\kc_{C\cap Y}(n)=\min\{d,\ \kc_Y(n)\}\ .\label{e6}
\end{equation}

At this moment we assume that $d>\sqrt{\deg(Y\cap C)}$. Taking
into account the choice of $d$ and property (\ref{e14}), and
looking at Figure \ref{f1}, we derive that
$$2d>\ord_1(Y)+1\quad\Longrightarrow\quad
d\ge\frac{\ord_1(Y)}{2}+1\ ,$$ and consequently
$$\deg(Y\cap C)\ge\frac{d(d+1)}{2}+\frac{(\ord_1(Y)+1-d)(\ord_1(Y)+2-d)}{2}
\ge\frac{(\ord_1(Y)+2)^2}{4}\ .$$ Hence, in view of $M_2(Y)\ge\deg
Y\ge\deg(Y\cap C)$,
$$\ord_1(Y)\le 2\sqrt{\deg(Y\cap
C)}-2<\sqrt{\frac{3}{2}M_2(Y)}+\frac{\deg Y}{\sqrt{3M_2(Y)/2}}-2\
.$$

{\bf Step 3}. From now on we assume that $Z$ is singular, and, in
the notation of Step 2
\begin{equation}
d\le\sqrt{\deg(Y\cap C)}\ .\label{e100}
\end{equation}

For the reader's convenience we start by proving a bound weaker
than (\ref{e7}),
\begin{equation}
\ord^{an}_1(Z) \le\sqrt{2M_2(Z)}+\frac{\deg Z}{\sqrt{2M_2(Z)}}-2\
.\label{e23}
\end{equation}
It will illustrate the main idea of the proof, which is based on
the argument used in the proof of Proposition \ref{p1}. Further
refinement up to inequality (\ref{e7}) is of technical nature and
consists of exploring particular steps in the proof of
(\ref{e23}).

\begin{remark}
The use of inequality (\ref{e23}) instead of (\ref{e7}) leads,
in fact, to similar estimates of orders of singular points
with a different constant factor.
\end{remark}

The choice of $d$ and properties (\ref{e16})-(\ref{e14}) of
Castelnuovo function immediately imply that
$$(\ord_1(Y)+1)d\le\deg(Y\cap C)+d(d-1)$$
\begin{equation}
\Longrightarrow\quad\ord^{an}_1(Z)\le\frac{\deg(Y\cap C)}{d}+d-2\
. \label{e5}
\end{equation}
We shall estimate $d=\deg C$ from below using the argument in the
proof of Proposition \ref{p1}.
The generality of $Y$ in $\Iso(Z)$ guarantees the fixed shape of
$\kc_Y$ when varying $Y$ in an open dense subset
of $\Iso(Z)$, as well as the
existence of a continuous family of curves $C$ of degree $d$,
satisfying (\ref{e6}) and having the same collection of
degrees of irreducible components. Let
$C=C_1^{l_1}...C_r^{l_r}$, where $C_1,...,C_r$ are reduced irreducible,
$\deg C_i=d_i$, $i=1,...,r$. Then, reasoning as
in the proof of Proposition
\ref{p1}, we obtain
$$d_i^2\ge\frac{1}{2}\sum_{q\in T(Y_{cl})}(\mt(C_i,q))^2,\quad i=1,...,r.$$
Taking generic elements $f_j\in I(Y_{w_j})$, $Y=Y_{w_1}\cup
...\cup Y_{w_p}$, $w_1,...,w_p\in\PP^2$, one obtains
$$\deg(C\cap Y)\le\sum_{j=1}^p(f_j\cdot C)_{w_j}
\le\sqrt{\sum_{q\in T(Y_{cl})}(\mt(Y,q))^2}\sqrt{\sum_{q\in T(Y_{cl})}
(\mt(C,q))^2}$$
$$\le\sqrt{2M_2(Y)}\cdot d\quad
\Longrightarrow\quad d\ge\frac{\deg(C\cap Y)}{\sqrt{2M_2(Y)}}\ .$$
In view of this bound and (\ref{e100}), inequality (\ref{e5})
implies (\ref{e23}).

{\bf Step 4}. To refine inequality (\ref{e23}) we need a
strengthened form of Lemma \ref{l1}.

\begin{lemma}\label{l2}
In the notation of Step 3, either $C_i$ is a straight line
containing two points of $T(Y_{cl})$, or
\begin{equation}
d_i^2\ge\frac{2}{3}\sum_{q\in T(Y_{cl})}(\mt(C_i,q))^2\ .\label{e25}
\end{equation}
\end{lemma}

The proof is found in section \ref{sec1} below.

We shall now prove (\ref{e7}) for schemes $Z$ such that,
for any two points $q_1,q_2\in T(Z_{cl})$,
\begin{equation}
\mt(Z_{cl},q_1)+\mt(Z_{cl},q_2)\le\sqrt{\frac{3}{2}M_2(Z)}\ .\label{e34}
\end{equation}

In the notation of Step 2, let $C=C'C''$,
$C'=C_1^{l_1}...C_s^{l_s}$, $C''=C_{s+1}^{l_{s+1}}...C_r^{l_r}$,
where $0\le s\le r$; any $C_i$, $i=1,...,s$, is a straight line
containing exactly two points of $T(Y_{cl})$, and
$C_{s+1},...,C_r$ are the other irreducible components of $C$.

Taking generic elements $f_j\in I(Y_{w_j})$, $Y=Y_{w_1}\cup
...\cup Y_{w_p}$, $w_1,...,w_p\in\PP^2$, one obtains
$$\deg(Y\cap C'')\le\sum_{j=1}^p(f_j\cdot C'')_{w_j}
\le\sqrt{\sum_{q\in T(Y_{cl})}(\mt(Y_{cl},q))^2}\sqrt{\sum_{q\in T(Y_{cl})}
(\mt(C'',q))^2}\ .$$
By Lemma \ref{l2},
$$\sqrt{\sum_{q\in T(Y_{cl})}(\mt(C'',q))^2}=
{\sum_{q\in T(Y_{cl})}\left(\sum_{i=s+1}^rl_i\cdot\mt(C_i,q)\right)^2}$$
$$\le\sum_{i=s+1}^rl_i\sqrt{\sum_{q\in T_i}(\mt(C_i,q))^2}
\le\sqrt{\frac{3}{2}}\sum_{i=s+1}^rl_id_i\ .$$
Hence
$$\deg(Y\cap C'')\le\sqrt{\frac{3}{2}M_2(Z)}\sum_{i=s+1}^rl_id_i\ .$$
On the other hand, by (\ref{e34}),
$$\deg(Y\cap C')\le\sum_{j=1}^p(f_j\cdot C')_{w_j}
=\sum_{i=1}^sl_i\left(\sum_{j=1}^p(f_j\cdot C_i)_{w_j}\right)$$
$$\le\sqrt{\frac{3}{2}M_2(Z)}\sum_{i=1}^sl_i\ .$$
So, it follows that
$$\deg(Y\cap C)\le\deg(Y\cap C')+\deg(Y\cap C'')
\le d\sqrt{\frac{3}{2}M_2(Z)}$$
$$\Longrightarrow\quad d\ge\frac{\deg(Y\cap C)}{\sqrt{3
M_2(Z)/2}}\ ,$$ which implies (\ref{e7}) by virtue of (\ref{e5})
and (\ref{e100}).

{\bf Step 5}. We shall complete the proof of Proposition \ref{p3}
by induction on $\deg Z$. The case of nonsingular $Z$ of degree
$\ne 2$ is the base of induction since (\ref{e22}) implies
(\ref{e7}) except for $\deg Z=2$.

Assume that $Z$ is singular. By the result of Step 4, one has to
consider only the case that $T(Z)$ has at least two vertices, and
there exist $q_1,q_2\in T(Z)$ such that
\begin{equation}
\mt(Z_{cl},q_1)+\mt(Z_{cl},q_2)>\sqrt{\frac{3}{2}M_2(Z)}\ .\label{e35}
\end{equation}
This means that, for a generic $Y\in\Iso(Z)$, there exists a straight line
$L$ containing two points of $T(Y)$ and satisfying
$$\deg(Y\cap L)>\sqrt{\frac{3}{2}M_2(Z)}=\sqrt{\frac{3}{2}M_2(Y)}\ .$$
In other words, $L$ contains two vertices $q_1,q_2$ of $T(Y_{cl})$
with multiplicities $m_i=\mt(Y_{cl},q_i)$, $i=1,2$, $m_1\ge 2$, and
$$\deg(Y\cap L)\ge m_1+m_2>\sqrt{\frac{3}{2}\sum_{q\in T(Y_{cl})}
(\mt(Y_{cl},q))^2}\ .$$

Put \begin{equation} d=\left[\sqrt{\frac{3}{2}M_2(Z)}+\frac{\deg
Z}{\sqrt{3M_2(Z)/2}}-2\right] \ .\label{e402} \end{equation}
First, notice that
\begin{equation}
d\ge m_1+m_2-1\ .\label{e36}
\end{equation}
Indeed,
$$\sqrt{\frac{3}{2}M_2(Z)}+\frac{\deg Z}{\sqrt{3M_2(Z)/2}}-2\ge
\sqrt{\frac{3}{2}M_2(Z)}+\frac{(m_1(m_1+1)+m_2(m_2+1))/2}{\sqrt{3M_2(Z)
/2}}-2$$
$$\ge\sqrt{\frac{3}{2}(m_1^2+m_2^2)}+\frac{m_1(m_1+1)+m_2(m_2+1)}{\sqrt{6
(m_1^2+m_2^2)}}-2\ge m_1+m_2-1\ ,$$ where the latter inequality
holds, because substituting $m_1$ and $m_2$ for $m=(m_1+m_2)/2$,
we diminish the left-hand side and obtain an inequality
$$\frac{4}{\sqrt{3}}m+\frac{1}{\sqrt{3}}-2\ge 2m-1\ ,$$
which holds true for $m\ge 3/2$.

The exact sequence of sheaves
$$0\to\kj_{Y:L/\PP^2}(d-1)\to\kj_{Y/\PP^2}(d)\to
\kj_{Y\cap L/L}(d)\to 0$$
induces the exact cohomology sequence
$$H^1(\kj_{Y:L/\PP^2}(d-1))\to H^1(\kj_{Y/\PP^2}(d))\to
H^1(\kj_{Y\cap L/L}(d))\ ,$$ where the last term vanishes due to
(\ref{e36}). Hence $H^1(\kj_{Y/\PP^2}(d))=0$, provided
$H^1(\kj_{Y:L/\PP^2}(d-1))=0$.

If $\deg(Y:L)\ne 2$, by the induction assumption, inequality
(\ref{e35}), and the formula $\deg(Y:L)=\deg Y-m_1-m_2$,
$$\ord_1(Y:L)\le
\sqrt{\frac{3}{2}M_2(Y:L)}+\frac{\deg Y:L}{\sqrt{3M_2(Y:L)/2}}-2$$
$$\le\sqrt{\frac{3}{2}M_2(Y)}+\frac{\deg Y-m_1-m_2}{\sqrt{3M_2(Y)/2}}-2
<\sqrt{\frac{3}{2}M_2(Y)}+\frac{\deg Y}{\sqrt{3M_2(Y)/2}}-3$$
$$\Longrightarrow\quad\ord_1(Y:L)\le d-1\ ,$$
and we are done.

If $\deg(Y:L)=2$, then, clearly, $H^1(\kj_{Y:L/\PP^2}(d-1))=0$
holds, provided, $d\ge 2$, which follows from (\ref{e402}), since
$m_1\ge 2$, $m_2\ge 1$, and consequently $\deg Z\ge 4$. \proofend

\subsection{Proof of Lemma \ref{l2}}\label{sec1}
Without loss of generality we assume that $T(Y_{cl})$ has at least
two vertices, $C_i$ contains at least two points of $T(Y_{cl})$,
and is not a straight line, containing exactly two points of
$T(Y_{cl})$.

Let $C_i$ contain exactly two points, say $q_1,q_2$ of
$T(Y_{cl})$, and is singular at $q_1,q_2$, then $d_i\ge 2$, and
(\ref{e25}) turns into the inequality $d_i^2\ge 4\ge 2/3(1+1)$.

Let $C_i$ contain exactly two points of $T(Y_{cl})$. If one of
these points $q$ is infinitely near to the other point $z$, and
$C_i$ is singular at $z$, i.e., $\mt(C_i,z)\ge 2$, then,
intersecting the curves $C_i$ and $\psi_1(C_i)$ (see notation in
the proof of Lemma \ref{l1}), we obtain by \cite{GS}, Theorem
2(1), or \cite{Xu}, Lemma 3,
$$d_i^2\ge(\mt(C_i,z))^2+(\mt(C_i,q))^2-\mt(C_i,q)\ge\frac{2}{3}
((\mt(C_i,z))^2+(\mt(C_i,q))^2)\ ,$$
where the latter inequality follows from
$$(\mt(C_i,z))^2+(\mt(C_i,q))^2-3\cdot\mt(C_i,q)
\ge 4+(\mt(C_i,q))^2-3\cdot\mt(C_i,q)>0\ .$$ If $C_i\cap
T(Y_{cl})$ consists of points $z_1\ne z_2\in\PP^2$, and $C_i$ is
singular at $z_1$, i.e., $\mt(C_i,z_1)\ge 2$, then, moving the
point $z_2$ and respectively the curve $C_i$, we similarly obtain
$$d_i^2\ge(\mt(C_i,z_1))^2+(\mt(C_i,z_2))^2-\mt(C_i,z_2)\ge\frac{2}{3}
((\mt(C_i,z_1))^2+(\mt(C_i,z_2))^2)\ .$$

The previous argument in general leads to (\ref{e25}), when
$C_i\cap Y$ is reducible. So, we assume that $C_i\cap Y$ is
irreducible, and $C_i$ contains at least three points of
$T(Y_{cl})$.

Let $k\ge 2$ in (\ref{e26}).
We shall refine the statement
of Lemma \ref{l1} up to (see the notation of
Lemma \ref{l1})
\begin{equation}
\sum_{q\in T_i\backslash T(\psi_k(Y))}
(\mt(C_i,q))^2\le\frac{1}{k+1}
\sum_{q\in T_i}(\mt(C_i,q))^2\ .\label{e27}
\end{equation}
Indeed, it would immediately follow from
\begin{equation}
\sum_{q\in T_i\cap T(\psi_{m+1}(Y))\backslash T(\psi_m(Y))}
(\mt(C_i,q))^2\le\sum_{q\in T_i\cap T(\psi_m(Y))\backslash
T(\psi_{m-1}(Y))}(\mt(C_i,q))^2\label{e29}
\end{equation}
for all $m\ge 1$. The vertices of $T(C_i)$, $T(Y)$ are encoded by
finite segments of the Puiseux expansion in a neighborhood of
$z=(0,0)$ $$C_i(x,y)=(1+O(x,y))\prod_{s=1}^n(y-\xi_s(x)),\quad
n=\mt(C_i,z),$$
$$f(Y)(x,y)=(1+O(x,y))\prod_{s=1}^{m_0}(y-\eta_s(x)), \quad
m_0=\mt(Y,z)\ ,$$ where $\xi_s(x)$, $\eta(x)$ are fractional power
series. For any $s=1,...,n$, there exists $p(s)$ such that
$$\ord(\xi_s(x)-\eta_{p(s)}(x))= \max_j\ord(\xi_s(x)-\eta_j(x))\ .$$
Then $$\sum_{q\in
T_i}(\mt(C_i,q))^2=\sum_{s=1}^n\sum_{t=1}^n\min\{\ord(\xi_t(x)-\xi_s(x),
\ord(\xi_t(x)-\eta_{p(s)})\}$$ $$=
\sum_{s=1}^n\sum_{t=1}^n\ord\left(\xi_t(x)-
\frac{\xi_s(x)+\eps_1\eta_{p(s)}(x)}{1+\eps_1}\right)\ ,$$ for a
generic number $\eps_1$. Similarly, $$\sum_{q\in T_i\cap
T(\psi_mY)}(\mt(C_i,q))^2=
\sum_{s=1}^n\sum_{t=1}^n\ord\left(\xi_t(x)-
\frac{\xi_s(x)+\eps_1\eta_{p(s)}(x)}{1+\eps_1}-\eps x^m\right)\
,$$ for generic numbers $\eps_1$, $\eps$. In view of the last
relation, the following series of simple inequalities yields
(\ref{e29}) and thereby (\ref{e27}): $$\ord\left(\xi_t(x)-
\frac{\xi_s(x)+\eps_1\eta_{p(s)}(x)}{1+\eps_1}-\eps
x^{m+1}\right)$$ $$\qquad\qquad-\ord\left(\xi_t(x)-
\frac{\xi_s(x)+\eps_1\eta_{p(s)}(x)}{1+\eps_1}-\eps x^m\right)$$
$$\le\ord\left(\xi_t(x)-
\frac{\xi_s(x)+\eps_1\eta_{p(s)}(x)}{1+\eps_1}-\eps x^m\right)$$
$$\qquad\qquad-\ord\left(\xi_t(x)-
\frac{\xi_s(x)+\eps_1\eta_{p(s)}(x)}{1+\eps_1}-\eps
x^{m-1}\right),$$ $s,t=1,...,n$.

In turn (\ref{e27}) and the assumption $k\ge 2$ yield
$$d_i^2\ge\sum_{q\in T_i\cap T(\psi_k(Y))}(\mt(C_i,q))^2
\ge\frac{k}{k+1}\sum_{q\in T_i}(\mt(C_i,q))^2
\ge\frac{2}{3}\sum_{q\in T_i}(\mt(C_i,q))^2\ .$$

The final case in the proof of Lemma \ref{l2} is that of $k=1$
being defined by (\ref{e26}). This yields, in particular, that
$C_i$ is singular at $z$. Indeed, $C_i$ is supposed to contain
three points of $T(Y_{cl})$, and if it were nonsingular, the third
point in $T(C_i)\cap T(Y_{cl})$ would have moved under
transformation $\psi_2$ which contradicts the assumption $k=1$.
Next, due to $T_i=T(C_i)\cap T(Y_{cl})\subset T(C_i)\cap
T(\psi_2(C_i))$, we obtain
\begin{equation}
\sum_{q\in T_i}(\mt(C_i,q))^2\le(C_i\cdot\psi_2(C_i))_z\
.\label{e30}
\end{equation}
Let $c(x,y)=0$ be
an equation of $C_i$ in an affine neighborhood of the
point $z=(0,0)$. Then
$$(C_i\cdot\psi_2(C_i))_z=\bigg(c(x,y)\cdot c(x,y+\eps x^2)
\bigg)_z$$
$$=\bigg(c(x,y)\cdot(c(x,y+\eps x^2)-c(x,y))\bigg)_z$$
$$\le\bigg(c(x,y)\cdot(x^2\frac{\partial c}{\partial y}(x,y))
\bigg)_z$$
\begin{equation}
=2\bigg(c(x,y)\cdot x\bigg)_z+\bigg(c(x,y)\cdot
\frac{\partial c}{\partial y}(x,y)\bigg)_z\ .\label{e31}
\end{equation}
On the other hand, for small $|\eps|$, in some neighborhood $U$ of
$z$, one has by \cite{GS}, Theorem 2(2),
$$d_i^2\ge(C_i\cdot\psi_1C_i)_U\ge\bigg(c(x,y)
\cdot(x\frac{\partial c}{\partial y}(x,y))\bigg)_z$$
\begin{equation}
=\bigg(c(x,y)\cdot x\bigg)_z+\bigg(c(x,y)\cdot
\frac{\partial c}{\partial y}(x,y)\bigg)_z\ .\label{e32}
\end{equation}
Here $\partial c/\partial y(0,0)=0$, since $C_i$ is singular at
$z$, and hence
\begin{equation}
\bigg(c(x,y)\cdot
\frac{\partial c}{\partial y}(x,y)\bigg)_z\ge
\bigg(c(x,y)\cdot x\bigg)_z\ .\label{e33}
\end{equation}
Combining (\ref{e30})-(\ref{e33}), one easily derives (\ref{e25}).

\section{Analytic and topological order of a singular point.
Curves with prescribed singularities}

\subsection{Plane curves with prescribed singularities}\label{sec6}

\begin{theorem}\label{t1}
Let $(C_i,z_i)$, $i=1,...,r$, be reduced plane curve germs with
isolated singular points $z_1,...,z_r$, respectively.

(1) If
\begin{equation}
\left[\sqrt{\frac{3}{2}\sum_{i=1}^rM_2(Z^s(C_i,z_i))}+
\frac{\sum_{i=1}^r\deg Z^s(C_i,z_i)}{\sqrt{3
\sum_{i=1}^rM_2(Z^s(C_i,z_i))/2}}\right]\le d+1\ ,\label{e39}
\end{equation}
then there exists an irreducible plane curve $C$ of degree $d$
having $r$ singular points topologically equivalent to
$(C_1,z_1)$, ..., $(C_r,z_r)$, respectively, as its only
singularities.

(2) If
\begin{equation}
\left[\sqrt{\frac{3}{2}\sum_{i=1}^rM_2(Z^a(C_i,z_i))}+
\frac{\sum_{i=1}^r\deg Z^a(C_i,z_i)}{\sqrt{3
\sum_{i=1}^rM_2(Z^a(C_i,z_i))/2}}\right]\le d+1\ ,\label{e43}
\end{equation}
then there exists an irreducible plane curve $C$ of degree $d$
having $r$ singular points analytically equivalent to $(C_1,z_1)$,
..., $(C_r,z_r)$, respectively, as its only singularities.

Furthermore, the germ at $C$ of the (topological or analytic)
equisingular stratum in the space of curves of degree $d$ is
T-smooth.
\end{theorem}

A particular case of one singularity is of special importance, since
it will be used in constructing curves with prescribed
singularities on arbitrary algebraic surfaces.

\begin{definition}\label{d1}
Given a reduced plane curve germ $(C,z)$, denote by $e^s(C,z)$
(resp., $e^a(C,z)$) the minimal degree $m$ of a plane curve $F$
having only one singular point $w$, which is topologically (resp.,
analytically) equivalent to $(C,z)$, and satisfying the condition
\begin{equation}
H^1(\kj_{Z/\PP^2}(m-1))=0\ ,\label{e70}
\end{equation}
where $Z=Z^{es}(F,w)$ (resp., $Z=Z^{ea}(F,w)$). We call the
parameters $e^s$ and $s^a$ the {\it topological} and {\it analytic}
order of a singular point.
\end{definition}

It should be noticed that $e^s$, $e^a$ introduced above differ
from similar singular point invariants used in
\cite{GLS,KT,Lo,Sh2}. The present notion corresponds to strong
transversality in the sense of \cite{Sh2}. More precisely,

\begin{lemma}\label{l13}
(1) Let $F$ be a plane curve as in Definition \ref{d1} and $L$ be
a straight line which does not pass through $w$. Then the germ at
$F$ of the family of curves of degree $m$ having in a neighborhood
of $w$ a singular point topologically (resp., analytically)
equivalent to $(C,z)$, is smooth of expected dimension, and
transversally intersects the linear system
$$\{G\in|\ko_{\PP^2}(d)|\ :\ G\cap L=F\cap L\}\ .$$

(2) Let $L\subset\PP^2$ be a straight line. Then the set of
$m$-tuples $\overline z$ of distinct points on $L$, such that
there is a curve $F$ of degree $m$ as in Definition \ref{d1},
satisfying $F\cap L=\overline z$, is Zariski open in $\Sym^m(L)$.
\end{lemma}

{\it Proof}.
The first statement is equivalent to
$$H^1(\kj_{Z\cup(F\cap L)/\PP^2}(m))=0$$
(see details in \cite{GL,GLS1,Sh}, where $Z$ is understood as in Definition \ref{d1}),
which follows from (\ref{e70}) and the exact sequence
$$0=H^1(\kj_{Z/\PP^2}(m-1))\to H^1(\kj_{Z\cup(F\cap L)/\PP^2}(m))
\to H^1(L,\kj_{F\cap L}(m))=0\ .$$

For the second statement take a curve $F$ as in Definition
\ref{d1}, which meets $L$ transversally, and the
germ $M$ at $F$ of the family of
curves of degree $m$ having in a neighborhood of $w$ a singular point
topologically (resp., analytically) equivalent to $(C,z)$.
Consider the map $G\in M\mapsto G\cap L\in\Sym^m(L)$.
Then the preceding statement of Lemma means that this map is
a submersion, and we are done.
\proofend

\begin{theorem}\label{t2}
For any reduced plane curve germ $(C,z)$,
\begin{equation}
e^s(C,z)\le\sqrt{\frac{3}{2}M_2(Z^s(C,z))}+ \frac{\deg
Z^s(C,z)}{\sqrt{3 M_2(Z^s(C,z))/2}}-1\ ,\label{e45}
\end{equation}
\begin{equation}
e^a(C,z)\le\sqrt{\frac{3}{2}M_2(Z^a(C,z))}+ \frac{\deg
Z^a(C,z)}{\sqrt{3 M_2(Z^a(C,z))/2}}-1\ .\label{e44}
\end{equation}
\end{theorem}

The hypotheses of Theorems \ref{t1}, \ref{t2} can be translated
into more familiar singularity invariants.

\begin{theorem}\label{t3}
Let $(C_i,z_i)$, $i=1,...,r$, and $(C,z)$ be plane curve germs
with isolated singular points $z_1,...,z_r$ and $z$, respectively.
Denote by $n$, $k$ and $t$ the number of nodes, cusps and
ordinary triple points, respectively,
among $(C_1,z_1)$, ..., $(C_r,z_r)$.

(1) If
\begin{equation}
6n+10k+\frac{169}{6}t+\frac{25}{3}u+\frac{27}{2}
 \sum_{(C_i,z_i)\ne A_1,A_2,D_4}\del(C_i,z_i)\le d^2-2d+3\ ,\label{e50}
 \end{equation}
 where $u$ is the number of points of type $A_{2m}$, $m\ge 2$,
 among $(C_i,z_i)$, $i=1,...,r$,
then there exists an irreducible curve of degree $d$ having $r$
singular points topologically equivalent to $(C_1,z_1)$, ...,
$(C_r,z_r)$, respectively, as its only singularities.

(2) If
 \begin{equation}
 6n+10k+\frac{169}{6}t+
 \sum_{(C_i,z_i)\ne A_1,A_2,D_4}\frac{(5\mu(C_i,z_i)+3\del(C_i,z_i)/2)^2}{3
\mu(C_i,z_i)+3\del(C_i,z_i)/2}\le d^2-2d+3\ ,\label{e57}
 \end{equation}
then there exists an irreducible curve of degree $d$ having $r$
singular points analytically equivalent to $C_1,z_1)$, ...,
$(C_r,z_r)$, respectively, as its only singularities.

(3) If $(C,z)$ is of type $A_m$, $m\ge 1$, then
 $$e^s(C,z)=e^a(C,z)\le 2[\sqrt{m+5}]\ ;$$
 if $(C,z)$ is of type $D_m$, $m\ge 4$, then
 $$e^s(C,z)=e^a(C,z)\le 2[\sqrt{m+7}]+1\ ;$$
 if $(C,z)$ is of type $E_m$, $m=6,7,8$, then
 $$e^s(C,z)=e^a(C,z)=\left[\frac{m+2}{2}\right]\ ;$$
 if $(C,z)$ is not simple, then
 $$e^s(C,z)\le\frac{9}{\sqrt{6}}\sqrt{\del(C,z)}-1\ ,$$
 $$e^a(C,z)\le\frac{5\mu(C,z)+3\del(C,z)/2}{\sqrt{3\mu(C,z)+3\del(C,z)/2}}-1
\le 3\sqrt{\mu(C,z)}-1\ .$$
\end{theorem}

\begin{remark}\label{r2}
We should like to point out that, for specific singularities,
invariants in the existence conditions can be reduced, for
example, a curve with $n$ nodes and $k$ cusps exists if $2n+4k\le
d^2+O(d)$ \cite{Sh1}, Theorem 4.1. We, however, have focused on
obtaining a universal existence condition rather than an
optimality of singularity invariants, though our results
substantially improve all previously known general existence
conditions \cite{GLS,Lo}. For example, since $\del\le2\mu/3$ for
singularities different from nodes, cusps and ordinary triple
points, (\ref{e57}) yields the following weaker, but rather
simpler sufficient existence condition for an irreducible plane
curve with one or many singularities prescribed up to analytic
equivalence:
$$e^a(C,z)\le 3\sqrt{\mu(C,z)}-2\ (\text{as}\ \mu(C,z)>4),
\quad\sum_{i=1}^r\mu(C_i,z_i)\le\frac{1}{9}(d^2-2d+3)\ ,$$
which is much better than the previously known
sufficient condition for the existence of an irreducible
curve with singularities prescribed only
up to topological equivalence (see\cite{Lo,GLS}):
$$e^s(C,z)<\min\left\{\sqrt{29\mu(C,z)}+\frac{9}{2},\
\sqrt{41.4\mu(C,z)}-2\right\}\ ,$$
$$\sum_{i=1}^r\mu(C_i,z_i)\le\frac{1}{46}(d+2)^2\ .$$
\end{remark}

\subsection{Proof of Theorem \ref{t1}}
Consider, first, the case of topological equivalence of singular
points. By Lemma \ref{l4}, without loss of generality, we can
suppose that $Z^s(C_i,z_i)$ is a generic element of
$\Def(Z^s(C_i,z_i))$, $i=1,...,r$.

Inequality (\ref{e39}) means by Proposition \ref{p3} that
\begin{equation}
H^1(\kj_{Z/\PP^2}(d-1))=0\ ,\label{e48}
\end{equation}
where $Z=Z^s(C_1,z_1)\cup...\cup Z^s(C_r,z_r)$. Introduce the
schemes
\begin{itemize}
\item $Z\cup\{z\}$, where $z\in\PP^2\backslash
\{z_1,...,z_r\}$;
\item $Z^{(i)}=Z^s(C_1,z_1)\cup...\cup Z^s(C_{i-1},z_{i-1})
\cup Z^s_1(C_i,z_i)\cup Z^s(C_{i+1},z_{i+1})\cup...\cup
Z^s(C_r,z_r)$, $i=1,...,r$.
\end{itemize}
We claim that
\begin{eqnarray}
&H^1(\kj_{Z\cup\{z\}/\PP^2}(d))=0,\quad z\in\PP^2\backslash Z\ ,
\label{e46}\\
&H^1(\kj_{Z^{(i)}/\PP^2}(d))=0,\quad i=1,...,r\ .\label{e47}
\end{eqnarray}
Indeed, let $L_0$ be a generic straight line through $z$, and $L_i$
be a generic
straight line through $z_i$. Then
$$\deg(L_0\cap Z^{(0)}(z))=\deg\{z\}=1\le d+1\quad
\Longrightarrow\quad H^1(L_0,\kj_{z}(d))=0\ .$$
In view of (\ref{e48}),
$$H^1(L_i,\kj_{L_i\cap Z}(d-1))=0\quad
\Longrightarrow\quad\deg(L_i\cap Z^s(C_i,z_i))\le d\ ;$$
hence
$$\deg(L_i\cap Z^s_1(C_i,z_i))
\le d+1\quad\Longrightarrow H^1(L_i,\kj_{Z^{(i)}}(d))=0\ .$$ Then,
using the last $h^1$-vanishing statement, the relation
$Z^s_1(C_i,z_i):L_i=Z^s(C_i,z_i)$ coming from Lemmas 2.14, 2.15
\cite{GLS}, and the exact sequences
\begin{eqnarray}
&0\to\kj_{Z/\PP^2}(d-1)\to\kj_{Z\cup\{z\}/\PP^2}(d)\to\kj_z(d)\to 0\ ,\nonumber\\
&0\to\kj_{Z/\PP^2}(d-1)\to\kj_{Z^{(i)}/\PP^2}(d)\to\kj_{L_i\cap
Z^s_1(C_i,z_i)/L_i} (d)\to 0\ ,\nonumber
\end{eqnarray}
one can easily derive (\ref{e46}), (\ref{e47}) from (\ref{e48}).

By Lemma \ref{l5}(2), the $h^1$-vanishing (\ref{e47}) implies
that, for any $i=1,...,r$, there exists a curve
$D_i\in|\kj_{Z/\PP^2}(d)|$ such that $(D_i,z_i)$ is topologically
equivalent to $(C_i,z_i)$, $1\le j\le r$. Again by Lemma
\ref{l5}(2), for a generic curve $D$ of the linear system
$\lam_1D_1+...+\lam_rD_r$, $\lam_1,...,\lam_r\in\C$, the germs
$(D,z_i)$ are topologically equivalent to $(C_i,z_i)$,
$i=1,...,r$, respectively. Furthermore, we can suppose that $D$ is
reduced. Let $w_1,...,w_m$ be all singular points of $D$ outside
$z_1,...,z_r$. By virtue of (\ref{e46}), there exist curves
$D'_1,...,D'_m\in|\kj_{Z/\PP^2}(d)|$ such that $w_j\not\in D'_j$,
$j=1,...,m$. By Bertini's theorem, the singular locus of a generic
curve $D'$ of the linear system $\lam
D+\Lam'_1D'_1+...+\lam'_mD'_m$, $\lam,\lam'_1,...,\lam'_m\in\C$,
is $\{z_1,...,z_r\}$, and $(D',z_i)$ is topologically equivalent
to $(C_i,z_i)$, $i=1,...,r$.

Finally, we show that $D'$ is irreducible. The above argument
shows that the linear system $|\kj_{Z/\PP^2}(d)|$ has no fixed
part. If all the curves in $|\kj_{Z/\PP^2}(d)|$ are reducible,
then by Bertini's theorem a generic curve
$D'\in|\kj_{Z/\PP^2}(d)|$ splits into irreducible components,
which all belong to a one-dimensional linear system. In
particular,
$$\dim|\kj_{Z/\PP^2}(d)|\le d\ .$$
On the other hand, by (\ref{e48}),
$$\dim|\kj_{Z/\PP^2}(d)|=\frac{d(d+3)}{2}-\deg Z\ .$$
Inequalities (\ref{e38}) and (\ref{e39}) yield
$$d+2>\sqrt{\frac{3}{2}M_2(Z)}+\frac{\deg Z}{\sqrt{3M_2(Z)/2}}$$
$$\ge\sqrt{\frac{3}{2}\deg Z}+\frac{\deg Z}{\sqrt{3\deg Z/2}}=
\frac{7}{\sqrt{6}}\sqrt{\deg Z}\ .$$
Hence
$$\dim|\kj_{Z/\PP^2}(d)|-d=\frac{d(d+1)}{2}-\deg Z$$
$$>\frac{13}{12}\deg Z-\frac{15}{2\sqrt{6}}\sqrt{\deg Z}+1\ge 0$$
as $\deg Z\ge 6$, and we are done, since the remaining cases of
one node or one cusp are obviously covered by (\ref{e39}).

The case of the analytic equivalence of singular points can be
treated in the same way, when using Lemmas \ref{l7}(2) and
\ref{l8} instead of
Lemma \ref{l5}(2).

Finally, we note that, by construction, $Z^s(D,z_i)=Z^s(C_i,z_i)$
in the first case of the Theorem, and $Z^a(D,z_i)=Z^a(C_i,z_i)$ in
the second case, $i=1,...,r$. In view of
$$I^s(D,z_i)\subset I^{es}(D,z_i),\quad I^a(D,z_i)\subset I^{ea}(D,z_i),
\quad i=1,...,r\ ,$$
(\ref{e48}) implies
$$H^1(\kj_{Z'/\PP^2}(d))=0\ ,$$
with the zero-dimensional scheme $Z'$ defined at points
$z_1,...,z_r$ by the ideals $I^{es}(D,z_i)$ or $I^{ea}(D,z_i)$,
respectively to the case considered. In turn the last
$h^1$-vanishing means the T-smoothness of the topological or
analytic equisingular stratum at $D$ in the space of curves of
degree $d$ (see \cite{GL,GLS,GLS1,Sh} for details).

\subsection{Proof of Theorem \ref{t3}}
The case $d\le 3$ is trivial, so we assume that $d\ge 4$.

(1) If there are no nodes and cusps among $(C_i,z_i)$,
$i=1,...,r$,
 then for an ordinary triple point $(C_i,z_i)$
$$\deg Z^s(C_i,z_i)=6,\quad M_2(Z^s(C_i,z_i))=9\ ,$$
for a point of type $A_{2m}$, $m\ge 2$, by (\ref{e38}) and Lemma
\ref{l10}, $$\deg Z^s(C_i,z_i)=3\del(C_i,z_i)+2,\quad
M_2(Z^s(C_i,z_i))=4\del(C_i,z_i)+2,$$
 and for the rest of the germs,
  $$\deg Z^s(C_i,z_i)\le 3\del(C_i,z_i),\quad M_2(Z^s(C_i,z_i))\le
  4\del(C_i,z_i)\ .$$
  Hence
$$\left(\sqrt{\frac{3}{2}\sum_{i=1}^rM_2(Z^s(C_i,z_i))}+
\frac{\sum_{i=1}^r\deg Z^s(C_i,z_i)}{\sqrt{3
\sum_{i=1}^rM_2(Z^s(C_i,z_i))/2}}\right)^2$$
\begin{equation}
\le\frac{(9\sum_{(C_i,z_i)\ne D_4}\del(C_i,z_i)+39t/2+5u)^2}{6
\sum_{(C_i,z_i)\ne D_4}\del(C_i,z_i)+27t/2+3u}
\le\frac{27}{2}\sum_{(C_i,z_i)\ne
D_4}\del(C_i,z_i)+\frac{169}{6}t+\frac{25}{3}u \ .\label{e54}
\end{equation}
 Thus, (\ref{e50}) implies (\ref{e39}), and we are done.

Assume that $n+k>0$. Put
 $$s=\max\{t\ge 2\ :\ \sum_{i=1}^{t-2}(d-i)+1\le 3n+5k\}\ .$$
 Observe that by (\ref{e50}), $s\le d-1$.
Inequalities (\ref{e50}), (\ref{e54}) also yield
 $$\sqrt{\frac{3}{2}\sum_{(C_i,z_i)\ne A_1,A_2}M_2(Z^s(C_i,z_i))}+
\frac{\sum_{(C_i,z_i)\ne A_1,A_2}\deg Z^s(C_i,z_i)}{\sqrt{3
\sum_{(C_i,z_i)\ne A_1,A_2}M_2(Z^s(C_i,z_i))/2}}$$
$$\le\sqrt{\frac{27}{2}
\sum_{(C_i,z_i)\ne
A_1,A_2,D_4}\del(C_i,z_i)+\frac{169}{6}t+\frac{25}{3}u}$$
 $$\le\sqrt{d^2-2d+3-6n-10k}\le\sqrt{d^2-2d+3-2\sum_{i=1}^{s-2}(d-i)-2}$$
$$=\sqrt{d^2-2(s-1)d+(s-1)(s-2)+1}\le d-s+1\ ;$$
hence, by Proposition \ref{p3}, $$H^1(\kj_{Z'/\PP^2}(d-s-1))=0\
,$$ where $Z'$ is the part of $Z$ without nodes and cusps. Then we
derive (\ref{e39}), and thereby the first statement of Theorem
\ref{t3}, from

\begin{lemma}\label{l12}
Let $L\subset\PP^2$ be a straight line,
$X,Y\subset\PP^2$ be zero-dimensional schemes such that
$L\cap Y=\emptyset$, $X\cap Y=\emptyset$, and $X=X^{(1)}\cup X^{(2)}
\cup X^{(3)}\cup X^{(4)}$, where $X^{(1)}$ is a union of schemes of degree
$1$ (i.e., locally defined by the maximal ideal), $X^{(2)}$ is a union
of schemes of degree $2$ (i.e., locally defined by ideals
like $\langle y,x^2\rangle$), $X^{(3)}$ is a union of schemes from
$\Def(Z^s(\text{node}))$ (i.e., locally defined by the square of the
maximal ideal), $Z^{(4)}$ is a union of schemes from
$\Def(Z^s(\text{cusp}))$ (i.e., locally defined by ideals like
$\langle y^2,yx^2,x^3\rangle$). Assume that
\begin{enumerate}
\item[(i)] $H^1(\kj_{Y/\PP^2}(d-s-1))=0$, where
$$s=\max\{t\ge 2\ :\ \sum_{i=1}^{t-2}(d-i)+1\le\deg X\}\quad\text{and}
\quad s\le d-1\ ;$$
\item[(ii)] $X^{(4)}\cap L=\emptyset$,  $\deg(X\cap L)\le d$;
\item[(iii)] the components of $X$ which do not meet $L$ are placed in
$\PP^2\backslash(Y\cup L)$ in a general position.
\end{enumerate}
Then
\begin{equation}
H^1(\kj_{X\cup Y/\PP^2}(d-1))=0\ .\label{e55}
\end{equation}
\end{lemma}

{\it Proof of Lemma \ref{l12}}. We perform induction on $s$, using
the so-called ``Horace method" \cite{H} (see also \cite{GLS} and
Step 4 in the proof of Proposition \ref{p3} above).

If $s=2$, then $\deg X\le d-1$. We specialize all the
components of $X$ on the line $L$ with maximal possible
intersection with $L$. First, we note that (\ref{e55}) for the
specialized scheme $X$ implies the same relation for the original $X$
due to the semicontinuity of cohomology. Second,
$$\deg((X:l)\cap L)\le\deg(X\cap L)\le\deg X\le d-1\ ,$$
and $X\subset L^2$. Then (\ref{e55}) follows from the two
exact sequences
$$0=H^1(\kj_{Y/\PP^2}(d-3))\to H^1(\kj_{(X:L)\cup Y/\PP^2}(d-2))
\to H^1(L,\kj_{(X:L)\cap L}(d-2))=0\ ,$$
\begin{equation}
H^1(\kj_{(X:L)\cup Y/\PP^2}(d-2))\to H^1(\kj_{X\cup Y/\PP^2}(d-1))
\to H^1(L,\kj_{X\cap L}(d-1))=0\ .\label{e56}
\end{equation}

Let $s\ge 3$ and $d-1\le\deg(X\cap L)\le d$. Then
$H^1(L,\kj_{X\cap L}(d-1))=0$, and (\ref{e56}) reduces (\ref{e55})
to $H^1(\kj_{(X:L)\cup Y/\PP^2}(d-2))=0$, which holds by the
induction assumption. Indeed, one can easily check that $X:L$ is
the union of zero-dimensional schemes of the same four kinds as
$X$, the intersection $(X:L)\cap L$ is the union of schemes of
degree $1$, and $\deg((X:L)\cap L)<\deg(X\cap L)\le d$. At last,
since
$$\sum_{i=1}^{s-1}(d-i)+1>\deg X\ ,$$
then
$$\sum_{i=1}^{s-2}(d-1-i)+1>\deg X-(d-1)\ge\deg(X:L)\ ,$$
so
$$s-1\ge\max\{t\ge 2\ :\ \sum_{i=1}^{t-2}(d-1-i)+1\le\deg X\}\ .$$

Let $s\ge 3$ and $\deg(X\cap L)\le d-2$. If all the components
of $X$ are specialized on $L$ then we complete the proof as
in the case $s=2$. If there are components of $X$ out of $L$,
we specialize some of them on $L$ keeping three rules:
\begin{itemize}
\item $\deg(X\cap L)\le d$,
\item if $\deg(X\cap L)\le d-2$ then a component $X'$ of degree
$2$ of $X$ is specialized on $L$ so that $\deg(X'\cap L)=2$
(with tangency),
\item if $\deg(X\cap L)\le d-3$ then a component $X'\in
\Def(\text{cusp})$ is specialized on $L$ so that
$\deg(X'\cap L)=3$ (with tangency).
\end{itemize}
If we end up with $d-1\le\deg(X\cap L)\le d$, then notice that the
components of $X:L$ meeting $L$ are of degree $\le 2$ or belong to
$\Def(\text{node})$, and $\deg((X:L)\cap L)\le\deg(X\cap L)-1\le
d-1$; this allows us to complete the proof as in the preceding
paragraph. If we end up with $\deg(X\cap L)\le d-2$, and all
components of $X$ specialized on $L$ so that the components from
$\Def(\text{cusp})$ are tangent to $L$, then $X\subset L^2$ and we
complete the proof as in the case $s=2$. \proofend

(2) If there are no nodes and cusps among $(C_i,z_i)$,
$i=1,...,r$,
 then by (\ref{e38}) and Lemma \ref{l10},
  $$\deg Z^a(C_i,z_i)\le 2\mu(C_i,z_i),\ M_2(Z^a(C_i,z_i))\le
  2\mu(C_i,z_i)+\del(C_i,z_i),\ (C_i,z_i)\ne D_4\ .$$
  Hence
$$\left(\sqrt{\frac{3}{2}\sum_{i=1}^rM_2(Z^a(C_i,z_i))}+
\frac{\sum_{i=1}^r\deg Z^a(C_i,z_i)}{\sqrt{3
\sum_{i=1}^rM_2(Z^a(C_i,z_i))/2}}\right)^2$$
 $$\le\frac{(\sum_{(C_i,z_i)\ne D_4}(5\mu(C_i,z_i)+3
\del(C_i,z_i)/2)+39t/2)^2}{\sum_{(C_i,z_i)\ne D_4}(3\mu(C_i,z_i)+3
\del(C_i,z_i)/2)+27t/2}$$
 $$\le\sum_{(C_i,z_i)\ne D_4}\frac{(5\mu(C_i)+3\del(C_i,z_i)/2)^2}{3
\mu(C_i,z_i)+3\del(C_i,z_i)/2}+\frac{169}{6}t\ .$$
 Thus, (\ref{e57}) implies (\ref{e43}), and we are done.

 If $n+k>0$ we prove the second statement of Theorem \ref{t3}
 in the same way as in the first part of the proof.

 (3) Estimates for $e^s(C,z)$, $e^a(C,z)$ (increased by $1$ in view of
 condition (\ref{e70})), where $(C,z)$ is of type $A_k,D_k$,
 are taken from
 \cite{Lo}. For singularities $E_k$, $k=6,7,8$, one has classically known curves
 of degree $d=4,4,5$, respectively. Relation (\ref{e70}) in these cases
 holds by \cite{DPW}, Theorem 1.1(2), under the condition $k<3d-3$.
The other estimates follow from (\ref{e38}), Lemma \ref{l10},
Theorem \ref{t2}, and the inequality $\del(C,z)\le
2/3\cdot\mu(C,z)$ for a non-simple singular point.

\subsection{Curves with prescribed singularities on algebraic
surfaces}\label{sec4} Let $\Sig$ be a smooth algebraic surface and
$D\subset\Sig$ be a divisor with $\dim|D|>0$. To obtain a
criterion for the existence of a curve $C\in|D|$ with prescribed
singularities, we combine \cite{Sh2}, Theorem 1, which basically
reduces the problem to $h^1$-vanishing for the ideal sheaf of a
zero-dimensional subscheme of $\Sig$ defined by $e^s$ or
$e^a$-powers of local maximal ideals, \cite{KT}, Theorem 2.1,
which provides a numerical sufficient condition for the above
$h^1$-vanishing, and Theorem \ref{t3}(3) with upper bounds to
$e^s$ and $e^a$.

\begin{theorem}\label{t4}
Let $\Sig$ be a smooth projective surface, $D$ a divisor on $\Sig$
with $D-K_{\Sig}$ nef, and $L\subset\Sig$ a very ample divisor.
Let $(C_1,z_1)$, ..., $(C_r,z_r)$ be reduced singular germs of
plane curves, among them $n$ nodes and $k$ cusps.

(1) If
\begin{equation}
18n+32k+27\sum_{\del(C_i,z_i)>1}\del(C_i,z_i)\le(D-K_\Sig-L)^2\
,\label{e60}
\end{equation}
\begin{equation}
(D-L-K_\Sig)L>\frac{9}{\sqrt{6}}\max_{1\le i\le r}\sqrt{\del(
C_i,z_i)}+1\ ,\label{e61}
\end{equation}
 and, for any irreducible curve $B$ with $B^2=0$ and
 $\dim|B|_a>0$,
\begin{equation}
\frac{9}{\sqrt{6}}\max_{1\le i\le r}\sqrt{\del(C_i,z_i)}<
 (D-K_\Sig-L)B+1\ ,\label{e62}
\end{equation}
 then there exists an irreducible curve $C\in|D|$ with $r$
 singular points topologically equivalent to $(C_1,z_1)$, ...,
 $(C_r,z_r)$, respectively, as its only singularities.

 (2) If
$$18n+32k+18\sum_{\mu(C_i,z_i)>2}\mu(C_i,z_i)\le(D-K_\Sig-L)^2\
,$$ $$(D-L-K_\Sig)L>3\max_{1\le i\le r}\sqrt{\mu(C_i,z_i)}+1\ ,$$
 and, for any irreducible curve $B$ with $B^2=0$ and
 $\dim|B|_a>0$,
 $$3\max_{1\le i\le r}\sqrt{\mu(C_i,z_i)}<
(D-K_\Sig-L)B+1\ ,$$
 then there exists an irreducible curve $C\in|D|$ with $r$
 singular points analytically equivalent to $(C_1,z_1)$, ...,
 $(C_r,z_r)$, respectively, as its only singularities.
\end{theorem}

Here $|B|_a$ means the family of curves algebraically equivalent
to $B$.

{\it Proof}. {\it Step 1}. Since, $e^s(\text{node})=2$,
$e^s(\text{cusp})=3$, and by Theorem \ref{t3}(3),
 $$e^s(C_i,z_i)+1\le\frac{9}{\sqrt{6}}\sqrt{\del(C_i,z_i)}\quad
 \text{as}\quad\del(C_i,z_i)>1\ ,$$
(\ref{e60}) and (\ref{e62}) imply
 $$2\sum_{i=1}^r(e^s(C_i,z_i)+1)^2\le(D-K_\Sig-L)^2\
 ,$$
 $$(D-K_\Sig-L)B>\max_{1\le i\le r}e^s(C_i,z_i)\ ,$$
 for curves $B$ as in (\ref{e62}). By \cite{KT}, Theorem 2.1,
this yields
 $$H^1(\Sig,\kj_Z(D-L))=0\ ,$$
 where $Z\subset \Sig$ is a zero-dimensional scheme concentrated
 at generic points $w_1,...,w_r\in\Sig$ and defined by the ideals
$(\mm_{w_i})^{e^s(C_i,z_i)}$, $i=1,...,r$.

{\it Step 2}. Let $w_0$ be a generic point in
$\Sig\backslash\{w_1,...,w_r\}$, and $w$ be any point in
$\Sig\backslash\{w_1,...,w_r,w\}$. Since
 $L$ is very ample, there is a nonsingular connected curve in
 $|L|$ (which we further denote by $L$ as well),
 which passes through $w,w_0$ and, may be, through one of $w_1,...,w_r$. In
 the exact sequence
 $$0=H^1(\Sig,\kj_Z(D-L))\to H^1(\Sig,\kj_{Z\cup\{w,w_0\}}(D))\to
 H^1(L,\kj_{(Z\cap L)\cup\{w,w_0\}}(D))\ ,$$
 the latter term vanishes, since by (\ref{e61}),
 $$\deg\kj_{(Z\cap L)\cup\{w,w_0\}}(D)\ge DL-2-
 \max_{1\le i\le r}e^s(C_i,z_i)>L^2+LK_\Sig=2g(L)-2\ .$$
 Hence
 $$H^1(\Sig,\kj_{Z\cup\{w,w_0\}}(D))=0\ .$$
 In particular, there exists a curve in $|\kj_{Z\cup\{w_0\}}(D)|$
 which does not pass through $w$. Since $w_0$ is generic, and $w$
 is any point outside
 $w_0,...,w_r$, by Bertini's theorem a generic curve $D_0\in
 |\kj_{Z\cup\{w_0\}}(D)|$ is nonsingular outside
 $w_1,...,w_r$, and this linear system has no fixed part.
 In addition, $D_0$ is irreducible. Indeed,
otherwise, by Bertini's theorem, $D_0$ and all close curves in
$|\kj_{Z\cup\{w_0\}}(D)|$ would split into variable components
which belong to the same one-dimensional algebraic family, but
this contradicts the fact that $D_0$ is nonsingular at the fixed
point $w_0$.

{\it Step 3}. For any $i=1,...r$,
define a zero-dimensional scheme $Z_i$ which coincides with $Z$ at
$w_j$, $1\le j\le r$, $j\ne i$, and is given by the ideal
$(\mm_{w_i})^{e^s(C_i,z_i)+1}$ at $w_i$. We claim that
 \begin{equation}
 H^1(\Sig,\kj_{Z_i}(D))=0\ .\label{e63}
 \end{equation}
 Indeed, there exists a nonsingular curve in $|L|$ (which we again denote by $L$
 for the sake of notation) which passes through $w_i$ and does not
contain any point $w_j$,
 $j\ne i$. Then in the exact sequence
 $$0=H^1(\Sig,\kj_Z(D-L))\to H^1(\Sig,\kj_{Z_i}(D))\to
 H^1(L,\kj_{Z_i\cap L}(D))\ ,$$
 the last term vanishes since by (\ref{e61}),
 $$DL-\deg(Z_i\cap L)\ge DL-e^s(C_i,z_i)>L^2+LK_\Sig=2g(L)-2\ .$$

{\it Step 4}. Using (\ref{e63}), for any $i=1,...,r$, we can find
a curve $D_i\in|\kj_{Z/\Sig}(D)|$ which has an ordinary singular
point of multiplicity $e^s(C_i,z_i)$ at $w_i$, and, in addition,
in some fixed local coordinates in a neighborhood of $w_i$, the
$e^s(C_i,z_i)$-jet at $w_i$ of a local equation of $D_i$ is a {\it
generic}\footnote{Here ``generic" means that the object considered
can be chosen arbitrarily in a Zariski open subset of the whole
space of objects.} $e^s(C_i,z_i)$-form. Then, a generic curve
$\widetilde D$ in the linear system spanned by $D_0,D_1,...,D_r$,
is irreducible, nonsingular outside $w_1,...,w_r$, has an ordinary
singular point of multiplicity $e^s(C_i,z_i)$ at $w_i$,
$i=1,...,r$, and, finally, the $e^s(C_i,z_i)$-jet at $w_i$ of a
local equation of $\widetilde D$ in some fixed coordinates in a
neighborhood of $w_i$ is a generic $e^s(C_i,z_i)$-form,
$i=1,...,r$.

{\it Step 5}. By Lemma \ref{l13}(2), for any $i=1,...,r$, there is
an affine curve $F_i$ of degree $e^s(C_i,z_i)$ with its only singular point
topologically equivalent to $(C_i,z_i)$, and such that the leading form
of the defining polynomial coincides with the $e^s(C_i,z_i)$-jet at
$w_i$ of a local equation of $\widetilde D$. Now we apply
\cite{Sh2}, Theorem 1, to deform $\widetilde D$ into the
required curve $C\in|D|$ with prescribed singularities.

Namely, in the assertion of \cite{Sh2}, Theorem 1,
$S=S_1=S_3=S_4=S_5=\emptyset$, $S_2=\{w_1,...,w_r\}$,
the affine curves $F_1,...,F_r$ serve as deformation models for
the singular points of $\widetilde D$, so that these models are
strongly transversal with respect to the topological
equivalence of singular points in view of (\ref{e70}).
Furthermore, the $h^1$-vanishing condition in \cite{Sh2},
Theorem 1, reads as
$$H^1(\Sig,\kj_{\widetilde Z}(D))=0\ ,$$
where $\widetilde Z=\bigcup_iZ^{es}(\widetilde D,w_i)$, which
immediately follows, say, from (\ref{e63}), since
$Z_i\supset\widetilde Z$ for any $i=1,...,r$.

{\it Step 6}. The second statement of Theorem \ref{t4} can be
proven in the same way, when taking into account that, by
Theorem \ref{t3}(3), $e^a(C_i,z_i)\le3\sqrt{\mu(C_i,z_i)}-1$
as $\mu(C_i,z_i)>2$.
\proofend

\section{Analytic order of a critical point}\label{sec3}

Let $f\in\hat\ko_{\C^2,0}$, $\mu(f)<\infty$. Denote by $e^a(f)$
the minimal degree of a polynomial $p\in\C[x,y]$ right equivalent
to $f$ at $0$, i.e., there is a local diffeomorphism
$\varphi:(\C^2,0)\to(\C^2,0)$ such that $p=f\circ\varphi$ in a
neighborhood of the origin.

\begin{theorem}\label{t5}
If $f$ is of type $A_m$, $m\ge 1$, then
 $$e^a(f)\le 2[\sqrt{m+5}]-1\ ;$$
 if $f$ is of type $D_m$, $m\ge 4$, then
 $$e^a(f)\le 2[\sqrt{m+7}]\ ;$$
 if $f$ is of type $E_m$, $m=6,7,8$, then
 $$e^a(f)=\left[\frac{m+2}{2}\right]\ ;$$
 if $f$ is not simple, then
\begin{equation}
e^a(f)<\frac{4}{\sqrt{3}}\sqrt{3\mu(f)-2\cdot\mt(f)+2}-1\
.\label{e72}
\end{equation}
\end{theorem}

{\it Proof}. Since the case of a simple germ $f$ coincides with
the case of a simple curve germ $\{f=0\}$, we take corresponding
bounds from Theorem \ref{t3}(3) (note only that here there is no
need to increase the estimates from \cite{Lo} by $1$, since
(\ref{e70}) is not required). Let $f$ be not simple. We claim that
 $$\ord_1^{an}(Z(f))<\frac{4}{\sqrt{3}}\sqrt{3\mu(f)-2\cdot\mt(f)+2}-1\
 .$$
 Indeed, by (\ref{e73}), Lemma \ref{l10}, and Proposition \ref{p3}
 $$\ord_1^{an}(Z_0(f))\le\frac{\dim Z_0(f)}{\sqrt{3M_2(Z_0(f))/2}}+
 \sqrt{\frac{3}{2}M_2(Z_0(f))}-2$$
 $$<\frac{4}{\sqrt{3}}\sqrt{\deg Z_0(f)}-2\le
\frac{4}{\sqrt{3}}\sqrt{3\mu(f)-2\cdot\mt(f)+2}-2\ .$$ Hence there
is a germ $g\in\hat\ko_{\C^2,0}$ right equivalent to $f$ such that
$$H^1(\ko{Z_0(g)/\PP^2}(m))=0,\quad
m=\left[\frac{4}{\sqrt{3}}\sqrt{3\mu(f)-2\cdot\mt(f)+2}\right]-2\
.$$ For a generic straight line $L$ intersecting $Z_0(g)$, we have
 $$Z(g):L=Z_0(g),\quad \deg(Z(g)\cap L)\le\mt Z_0(g)+1\le\deg
 Z_0(g)\le m\ ;$$
 hence
 $$0=H^1(\ko_{Z_0(g)/\PP^2}(m))\to H^1(\ko_{Z(g)/\PP^2}(m+1))
 \to H^1(L,\ko_{Z(g)\cap L}(m+1))=0\ ,$$
 which yields
 $$H^1(\ko_{Z(g)/\PP^2}(m+1))=0\ .$$
 Then the surjectivity of the morphism
 $$H^0(\ko_{\PP^2}(m+1))\to H^0(\ko_{Z(g)})=\hat\ko_{\C^2,0}/I(g)\
 ,$$
 gives us a polynomial $p\in\C[x,y]$ of degree $\le m+1$
 such that $p-g\in I(g)$, and we are done by Lemma \ref{l11}(2).

\section{Higher-dimensional case: example}\label{sec7}

Let $f:(\C^n)\to(\C,0)$, $n\ge 3$, be a germ of a holomorphic
function with an isolated critical point, i.e., $\mu(f)<\infty$.
As for $n=2$, we would like to estimate the minimal degree $d(f)$
of a polynomial $p\in\C[x_1,...,x_n]$ right equivalent to $f$ at
the origin. The classical bounds are
 $$\sqrt[n]{\mu(f)}+1\le d(f)\le\mu(f)+1\ .$$
 Analogously to the two-dimensional case we state

 \begin{conjecture}\label{c1}
There exists a sequence of positive numbers $a_n$, $n\ge 3$, such
that
 $$d(f)\le a_n\sqrt[n]{\mu(f)}$$
 for any germ $f:(\C^n)\to(\C,0)$ with $\mu(f)<\infty$.
 \end{conjecture}

This conjecture is elementary for $n=1$, and follows from Theorem
\ref{t5} for $n=2$. It, in fact, implies similar bounds
for isolated singular points of hypersurfaces in $\PP^n$, and
sufficient existence conditions for hypersurfaces of a given degree with
prescribed isolated singularities.

To support Conjecture \ref{c1}
we prove it in any
dimension for the case of critical points of type
$A_k$, $k\ge 1$.

\begin{theorem}\label{t6}
There is a sequence of positive numbers $a_n$, $n\ge 1$, such that
\begin{equation}
d(f)\le a_n\sqrt[n]{k}\label{e80}
\end{equation}
 for any germ $f:(\C^n,0)\to(\C,0)$, $n\ge 1$, of type $A_k$,
 $k\ge 1$.
\end{theorem}

\begin{remark}
For $n=2$ one can produce explicit formulas for polynomials of
degree $\le a_2\sqrt{k}$ with a critical point of type $A_k$, like
the classically known polynomial $(y-x^m)^2+y^{2m}$ of degree $2m$
with the critical point of type $A_{2m^2-1}$ at the origin (see
more examples in \cite{GN,Lo}). We do not know similar formulas for
$n\ge 3$, and provide an existence proof in the spirit of
preceding sections.
\end{remark}

{\it Proof}. For a germ $f:(\C^n,0)\to(\C,0)$ of type $A_k$, $k\ge
1$, introduce the zero-dimensional schemes $Z_0(f)$ and $Z(f)$
defined at the origin by the ideals
 $$I_0(f)=\{g\in\hat\ko_{\C^n,0}\ :\ g,g_{x_1},...,g_{x_n}\in
 \langle f_{x_1},...,f_{x_n}\rangle\},\quad I(f)=\mm_0\cdot I_0(f)\
 ,$$
 (cf. section \ref{sec2}). Similarly to the proof of Theorem
 \ref{t5} we conclude that
 $$d(f)\le\ord_1^{an}(Z(f))\le\ord_1^{an}(Z_0(f))+1\ .$$
We shall show that there is a sequence of positive numbers $b_n$,
$n\ge 1$, such that
 \begin{equation}
 \ord_1^{an}(Z_0(f))\le b_n\sqrt[n]{k}\ .\label{e81}
 \end{equation}

 In the sequel we shall use an auxiliary statement.

 \begin{lemma}\label{l14}
 In the given notation:

(1) After a suitable linear coordinate change, the hypersurface
germs $\{f_{x_2}=0\}$, ..., $\{f_{x_n}=0\}$ are nonsingular and
intersect transversally along a nonsingular curve
$$\gam(t)=(x_1(t),...,x_n(t)),\quad x_1(t)=t,\ x_i(t)=O(t^2),\
i=2,...,n\ .$$ The intersection multiplicity of $\gam$ and
$\{f_{x_1}=0\}$ at the origin is
$$(\gam\cdot f_{x_1})_0=\mu(f)=k\ .$$

(2) Introduce an ascending sequence of zero-dimensional schemes
$Z_m$, $m\ge 1$, defined at the origin by the ideals
$$I_m=\{g\in\ko_{\C^n,0}\ :\ (\gam\cdot g)_0,(\gam\cdot
g_{x_1})_0,..., (\gam\cdot g_{x_n})_0\ge m\},\quad m\ge 1\ .$$ Let
$L\subset\PP^n$ be a hyperplane such that
$L\cap\C^n=\{x_n=0\}$ and $2\le(\gam\cdot
L)_0=m_0\le k$. Assume that the function germ
$f|_L:(L,0)\to(\C,0)$ has an isolated critical point. Then
\begin{enumerate}
\item[(i)] $Z_k=Z_0(f)$,
\item[(ii)] $Z_m\cap L\subset Z_0(f|_L)$ for all $m\le k$,
%\item[(iii)] $\deg(Z_m\cap L)\ge m_0$ if $m\ge m_0$,
\item[(iii)] $Z_m:L^2\subset Z_{m-m_0}$ if $m>m_0$, and
$Z_m\subset L^2$ if $m\le m_0$.
\end{enumerate}
 \end{lemma}

 {\it Proof of Lemma \ref{l14}}.
For the first claim we observe that for suitable coordinates
$f(x_1,...,x_n)=x_2^2+...+x_n^2+\text{h.o.t.}$

Statement 2(i) becomes trivial when passing to local coordinates
$y_1,...,y_n$ such that $f=y_1^{k+1}+y_2^2+...+y_n^2$. In
statement 2(ii) it is enough to check the case $m=k$, which
reduces to the evident implication
 $$g\in\langle
 f_{x_1},...,f_{x_n},x_n\rangle\quad\Longrightarrow\quad
 g|_L\in\langle (f_{x_1})|_L,...,(f_{x_{n-1}})|_L\rangle\ .$$
 Statement 2(iii) is straightforward. \proofend

 Now we are going to prove (\ref{e81}) by induction on $n$; the
 case $n=1$ is trivial, the case of $n=2$ is covered by Lemma
 \ref{l10} and Proposition \ref{p3}. Assume that $n\ge 3$ and that
 (\ref{e81}) is proven for $n-1$. Let $(s-1)^n<k\le s^n$. Put
 $m_0=s^{n-1}$. By the induction assumption there is a function
germ $h=h(x_1,...,x_{n-1}): (L,0)\to(\C,0)$ with an isolated
critical point of type $A_{2m_0-1}$ such that
\begin{equation}
H^1(L,\kj_{Z_0(h)}(p))=0,\quad p\ge[b_{n-1}2^{1/(n-1)}s]\ .\label{e82}
\end{equation}
By a linear coordinate change in $L$ we can turn $h$ into
$$h(x_1,...,x_{n-1})=x_2^2+...+x_{n-1}^2+\text{h.o.t.}$$
In particular, the system
$$h_{x_2}=...=h_{x_{n-1}}=0$$
has a solution $x_2=x_2(x_1)$, ..., $x_{n-1}=
x_{n-1}(x_1)$ in a neighborhood of $0$, and
$$\varphi(x_1)=h_{x_1}(x_1,x_2(x_1),...,x_{n-1}(x_1))=
\alp x_1^{2m_0-1}+O(x_1^{2m_0}),\quad\alp\ne 0\ .$$
Define a function germ $f:(\C^n,0)\to(\C,0)$ by
$$f(x_1,...,x_n)=h(x_1,...,x_{n-1})+x_n^2-2x_n\psi(x_1)$$
with certain $\psi:(\C,0)\to(\C,0)$. Here the system
$$f_{x_2}=...=f_{x_n}=0$$
defines a curve $\gam:(\C,0)\to(\C^n,0)$
$$x_1=t,\ x_2=x_2(t),\ ...,\ x_{n-1}=x_{n-1}(t),\ x_n=
\psi(t)\ ,$$
and we impose the condition $(\gam\cdot f_{x_1})_0=k$, saying that
$f$ has an isolated critical point of type $A_k$.
To satisfy the last condition,
we choose $\psi$ so that
$$f_{x_1}(t,x_2(t),...,x_{n-1}(t),\psi(t))=t^k\ ,$$
which in view of $f_{x_1}=h_{x_1}-2x_n\psi'(x_1)$ reads as
$$\varphi(t)-2\psi(t)\psi'(t)=t^k\ ,$$
and is solved with respect to $\psi$ as follows
$$\psi^2(t)=\int_0^t(\varphi(t)-t^k)dt=\frac{\alp}{2m_0}t^{2m_0}
+O(t^{2m_0+1})-\frac{t^{k+1}}{k+1}$$
$$\Longrightarrow\quad
\psi(t)=\sqrt{\frac{\alp}{2m_0}}t^{m_0}+O(t^{m_0+1})\ .$$ Observe
also that the last formula means that $(\gam\cdot L)_0=m_0$, and
thus that Lemma \ref{l14} applies.

We complete the proof by establishing
$$H^1(\kj_{Z_0(f)/\PP^n}(s_0+2s))=0,\quad s_0=[b_{n-1}2^{1/(n-1)}s]\ .$$
Indeed, consider the exact sequences $\ke_i$, $i=1,...,2s$,
$$\ke_i\ :\qquad H^1(\kj_{Z_0(f):L^i/\PP^n}(s_0+2s-i))\to
H^1(\kj_{Z_0(f):L^{i-1}/\PP^n}(s_0+2s-i+1))$$
$$\to H^1(L,\kj_{(Z_0(f):L^{i-1})\cap L}(s_0+2s-i+1))\ .$$
By Lemma \ref{l14}(2)(ii)
$$(Z_0(f):L^{i-1})\cap L\quad\subset\quad Z_0(f)\cap L
\quad\subset\quad Z_0(h)\ ;$$ hence the last term in all the
sequences vanishes by (\ref{e82}). Observing that
$Z_0(f):L^{2s}=\emptyset$ by Lemma \ref{l14}(2)(iii) and $k/m_0\le
s$, we conclude that all terms in the sequence $\ke_{2s}$ vanish,
and then moving inductively by decreasing $i$, we obtain that the
middle term in $\ke_1$, which appears in (\ref{e82}), vanishes.

\end{document}